\documentclass[11pt,myheadings]{article}
\newtheorem{thm}{Theorem}[section]
\newtheorem{prop}{Proposition}[section]
\newtheorem{defn}{Definition}[section]
\newtheorem{lemma}{Lemma}[section]
\newtheorem{cor}{Corollary}[section]
\newtheorem{rem}{Remark}[section]
\newtheorem{exam}{Example}[section]
\usepackage{amsmath}
\usepackage{amssymb}

\begin{document}
\title{\bf Factorization from a poset-theoretic view I
\thanks {Keywords and phases: $B$-ideals, divisorial ideals,
a Galois connection.}}
\author{\small Zike Deng }
\date{}
\maketitle \noindent{\bf Abstract:} We introduce B-ideals and based
on them establish several necessary and sufficient conditions for an
element of  a monoid to be decomposed into a least common multiple
of infinite or a finite number of powers of prime factors. Besides,
we introduce a sort of Galois connection relating them to divisorial
ideals.

\section{Introduction}
The objective of this treatise, which consists of two parts, is to
investigate factorization in a monoid by means of the relation of
divisibility mainly, while the multiplication plays a secondary
roll. Or rather, what relates to factorization is
something essentially poset-theoretic, which through connection of
order with operation transforms into an algebraic result. We are
partly inspired in ideas by [1,5] and in techniques by [2,6] so as
to do this. We introduce B-ideals, which connect decomposition with
complete distributivity [3], as the tool and establish a sort of
Galois connection relating them to divisorial integral ideals [1].

As early as in the thirties of last century it was recognized that
factorization in itself referred to multiplication and dispensed
with addition,however,factorization in a monoid is a problem which
has not yet been solved up to now because no ideal is available.Here
there is divisibility only that remains.Based on it we reduce the
original problem to the poset-theoretic one and solve it.Then we
transform the results back into the algebraic ones.In contrast,we
view an ideal as an element of a monoid and so its decomposition
also becomes factorization.B-ideal, in case of a domain,relates to a
divisorial integral ideal through a Galois connection and hence it
is in a sense a generalization of the latter to a monoid.B-ideal
also relates to the dual of a filtre in topology.$J_b$ is the
analogue of the dual of open neighborhood base, of the kernel of a
valuation restricted to its ring in the poset-theoretic setting
respectively. B-ideals can avoid the same set of factors determining
different ideals and different sets of factors determining the same
divisorial ideal. That is why they work in factorization.

We study arbitrary decomposition, i.e., a least common multiple of
infinite or a finite number of powers of prime factors, and lay
emphasis on a unique factorization domain (taking as a monoid the
collection of equivalence classes determined by the preorder of
divisibility) and a Krull domain (taking as a monoid the collection
of integral divisors).

We establish topological representation [4] to which the
decomposition leads directly, and introduce poset-theoretic
constructions such as subposets of the first kind (of the second
kind), internal or external (direct) products and their mutual
relations, which are not only concerned with order representation
but also connect the decomposition with structure problems of the
monoid.

From a poset-theoretic view arbitrary decomposition in terms of
irreducible elements is connected with topology,while its finiteness
with algebra.Order,so to speak,bridges the gap between the both.
Decomposition can be rephrased as follows.Does there exist a
topology such that G is the set of closed subsets with B the set of
point-closures? In finite case of order construction a product
corresponds algebraically to the same of a monoid as a semimodule
and a subposet of the first kind, of the second kind to a
projective,injective semimodule respectively. A factorial monoid G
is a free semimodule on positive integers with primes as a base, a
closed-set lattice with powers of primes as point-closures and a
completely distributive complete lattice.

Through this treatment we recognize that an operation generates both
a monoid and a preorder, and this preorder is related closely to the
structure of the monoid through poset-theoretic constructions and
under a certain condition generates a topology which associates
itself with the monoid intrinsically.

This paper is the first part, which is concerned with arbitrary
decomposition. Section 2 introduces B-ideals and by using them
treats the characterizations of decomposition of an element of a
monoid into an arbitrary join of powers of prime factors in terms of
the properties of B-ideals (Theorem 2.1). Besides, uniqueness of
decomposition (Proposition 2.6) and its topological representation
(Proposition 2.8) are also studied. In section 3, we establish a
Galois connection (Proposition 3.1) which relates integral divisors
to B-ideals so as to transform decomposition problem of the former
into that of the latter. At last we introduce poset-theoretic
constructions such as subposets of the first kind or of the second
kind, internal or external products(Examples 3.2-3.7).

\section{Arbitrary decomposition}
Throughout this treatise $R$ will be an integral domain, i.e., a
commutative ring with identity $1\neq 0$ and without zero divisors.
Put $R^\ast=R\setminus \{0\}$. The relation of divisibility $x|y$
(equivalent to $y=zx$ for some $z\in R^\ast$) makes $R$($R^\ast$)
into a preordered monoid, i.e. ,the preorder $|$ is compatible with
the multiplication, and $G=R/U$ ($G^\ast=R^\ast/U$), where $U$ is
the set of all units of $R$, is an ordered monoid. Denote an element
of $G$ by $[x]$, the equivalence class of $x\in R$ and we have
$G^\ast=G\setminus \{[0]\}$. Put $B=\{[p]^n\mid p\in R$ is prime and
$n=1,2,\cdots\}$ and now we will consider decomposition of an
element of $G^\ast$ into an arbitrary join of elements of $B$. By
join (meet) we mean supremum (infinimum), or least common multiple
(greatest common divisor) if the order is relation of divisibility.
Denote it by $\vee$ ($\wedge$). We will regard $G^\ast$ as a poset
only and for brevity an element of $G^\ast$ will be denoted by $a,
b, c$. For $a\in G^\ast$, $A\subseteq G^\ast$ denoted by $\downarrow
a$, $\downarrow A$ the set $\{c\in G^\ast\mid c\leqslant a\}$, $\cup
\{\downarrow a\mid a\in A\}$ respectively. $A$ with $\downarrow A=A$
is called a lower set. $2^{(G^\ast)}$ refers to the collection of
all subsets of $G^\ast$ with joins existing.

\begin{defn} $J\subseteq G^\ast$ with $J\neq \varnothing$ is called
a $B$-ideal if $\downarrow J=J$ and $2^{(G^\ast)}\ni \downarrow
a\cap B\subseteq J$ implies $\vee \downarrow a\cap B\in J$ for each
$a\in G^\ast$.
\end{defn}

Denote by $M$ the collection of all $B$-ideals. Since $G^\ast\in M$,
we might as well regard $G^\ast$ as $\downarrow [0]$ and hence for
each $a\in G$, $\downarrow a$ is a $B$-ideal and is called a
principal $B$-ideal. Note that $\downarrow [1]$ is the least element
of $M$ and $\downarrow [0]=G^\ast$ the greatest element of $M$.

Put $\uparrow a=\{c\mid G^\ast\ni c\geqslant a\}$ and $\uparrow
A=\cup \{\uparrow a\mid a\in A\}$ for $a\in G^\ast, A\subseteq
G^\ast$.

\begin{defn}
{\rm(1)} $A\subseteq G^\ast$ is called a $B$-set if $\uparrow A=A$
and for each $a\in A$ there exists $b\in B\cap A$ such that
$b\leqslant a$.

{\rm(2)} $A\subseteq G^\ast$ is called a $B$-filter if $A$ is a
$B$-set satisfying that for all $a,b\in A$ with $a\wedge b$
existing, $a\wedge b\in A$.

{\rm(3)}The $B$-ideal $J$ is said to be prime if for all $a,b\in
G^\ast$ with $a\wedge b$ existing,$a\wedge b\in J$ implies $a\in J$
or $b\in J$.
\end{defn}

Assume that $P$ is a poset. In case of $\vee A$ existing for any
nonempty finite subset $A\subseteq P$, $a\in P$ is said to be
(strongly) $\vee$-irreducible if $a=b\vee c$ ($a\leqslant b\vee c$)
implies $a=b$ or $a=c$ ($a\leqslant b$ or $a\leqslant c$) for all
$a, b, c\in P$. In case of $\vee A$ existing for any subset
$A\subseteq P$, $a\in P$ is said to be (strongly) completely
$\vee$-irreducible if $a=\vee A$ ($a\leqslant \vee A$) implies $a=b$
($a\leqslant b$) for some $b\in A$. Using $\geqslant, \wedge$
instead of $\leqslant, \vee$ respectively in above-mentioned
definitions we obtain corresponding definitions of being (strongly)
$\wedge$-irreducible and (strongly) completely $\wedge$-irreducible.

A poset $P$ with $\vee A$ ($\wedge A$) existing for any nonempty
finite subset $A\subseteq P$ is called a $\vee$-semilattice (a
$\wedge$-semilattice). A poset $P$ is called a lattice if it is both
a $\vee$-semilattice and a $\wedge$-semilattice. A poset $P$ with
$\vee A$, $\wedge A$ existing for any subset $A\subseteq P$ is
called a complete lattice. A lattice $P$ is said to be distributive
if $a\wedge (b\vee c)=(a\wedge b)\vee (a\wedge c)$ for all $a, b,
c\in P$. A complete lattice $P$ is said to be completely
distributive if $\wedge_{t\in T}\vee A_t=\vee_{f\in \prod A_t}\wedge
f(T)$ for any $\{A_t\mid t\in T\}\subseteq 2^P$.

\begin{prop}
{\rm(1)} If $G^\ast\setminus J$ is a $B$-set $($a $B$-filter$)$,
then $J$ is a $B$-ideal $($a prime $B$-ideal$)$.

{\rm(2)} If $J\in M$ is strongly $\wedge$-irreducible, then $J$ is
prime.
\end{prop}
{\bf Proof.} (1) Let $c\leqslant a\in J$, then $a\notin
G^\ast\setminus J$, a fortiori, $c\notin G^\ast\setminus J$ and
hence $c\in J$. Assume that $2^{(G^\ast)}\ni \downarrow a\cap
B\subseteq J$. If $\vee \downarrow a\cap B\notin J$, then $\vee
\downarrow a\cap B\in G^\ast \setminus J$ and so $b\leqslant \vee
\downarrow a\cap B$ for some $b\in (G^\ast \setminus J)\cap B$,
whence $b\leqslant a$ because $\vee \downarrow a\cap B\leqslant a$.
Thus $b\in \downarrow a\cap B\subseteq J$, a contradiction to $b\in
G^\ast \setminus J$. It follows that $\vee \downarrow a\cap B\in J$
and hence $J$ is a $B$-ideal. On the other hand, suppose $G^\ast
\setminus J$ is a $B$-filter. We have already known from above that
$J$ is a $B$-ideal. Let $a\wedge b\in J$ for $a,b\in G^\ast$ with
$a\wedge b$ existing. If $a,b\notin J$, then $a,b\in G^\ast
\setminus J$ and so $a\wedge b\in G^\ast \setminus J$, a
contradiction. Hence $J$ is prime.

(2) Let $a\wedge b$ exist and $a\wedge b\in J$. Then $\downarrow
a\cap \downarrow b=\downarrow (a\wedge b)\subseteq J$ and so
$\downarrow a\subseteq J$ or $\downarrow b\subseteq J$, whence $a\in
J$ or $b\in J$.\ \ $\square$

\begin{prop}
$M$ ordered by inclusion is a complete lattice.
\end{prop}
{\bf Proof.} Note first that $M$ has the least element $\downarrow
[1]$ and the greatest element $\downarrow [0]= G^\ast$. Let
$\{J_t\mid t\in T\}$ with $T\neq \varnothing$, then $\cap J_t\neq
\varnothing $ because $[1]\in J_t$ for each $t\in T$. It is clear
that $\cap J_t\in M$. \ \ $\square$

~

Evidently, $J=\vee_{a\in J}\downarrow a$ for each $J\in M$. For
$A\subseteq G^\ast$ put $J(A)=\cap\{J\mid A\subseteq J\in M\}$,
which is called the $B$-ideal generated by $A$. For any $\{J_t\mid
t\in T\}\subseteq M$ we have $\vee J_t=J(\cup J_t)$.

For $b\in B$ define $J_b=\{a\mid b\nleqslant a\in G^\ast\}$ and
denote by $\sum_1$ the collection of all such $J_b$'s. It is trivial
that $J_b$ is the greatest in all $B$-ideals missing $b$.

\begin{prop}
$J_b\in M$
\end{prop}
{\bf Proof.} Let $c\leqslant a\in J_b$, then $b\nleqslant a$, a
fortiori, $b\nleqslant c$ and hence $c\in J_b$. Besides, assume that
$2^{(G^\ast)}\ni \downarrow a\cap B\subseteq J_b$. If $\vee
\downarrow a\cap B\notin J_b$, then $b\leqslant \vee \downarrow
a\cap B\leqslant a$, whence $b\in \downarrow a\cap B$ and so $b\in
J_b$, a contradiction. Thus $\vee \downarrow a\cap B\in J_b$. It
follows that $J_b\in M$.\ \ $\square$

It will be readily verified that for each $b\in B, G^\ast \setminus
J_b=\uparrow b$ is a $B$-filter and by Proposition 2.1(1) $J_b$ is a
prime $B$-ideal.

\begin{prop}
$J_b$ is strongly completely $\wedge$-irreducible.
\end{prop}
{\bf Proof.} Assume that $J_b\supseteq \cap J_t$ for any $\{J_t\mid
t\in T\}\subseteq M$. Since $b\notin J_b$, we have $b\notin \cap
J_t$ and hence $b\notin J_{t_0}$ for some $t_0\in T$. Let $a\in
J_{t_0}$, then $b\nleqslant a$ because otherwise $b\in J_{t_0}$,
hence $a\in J_b$. Thus $J_{t_0}\subseteq J_b$.\ \ $\square$

Put $\vartriangle_a=\{J\mid J\in M$ and $a\notin J\}$.

\begin{defn}
{\rm(1)} By the condition $D_1$ we mean that $a=\vee \downarrow
a\cap B$ for each $a\in G^\ast$.

{\rm(2)} By the condition $D_2$ we mean that for all $a\in G^\ast,
J\in M, J\in \vartriangle_a$ implies $J\in \vartriangle_b$ and
$J_b\in \vartriangle_a$ for some $b\in B$.

{\rm(3)} By the condition $D_3$ we mean that $J=\cap_{b\in
B\setminus J}J_b$ for each $J\in M$.

{\rm(4)} By the condition $D_4$ we mean that $J=\vee_{b\in J\cap
B}\downarrow b$ for each $J\in M$.

{\rm(5)} By the condition $D_5$ we mean that every nonzero element
of $R$ can be written as a least common multiple of infinite or a
finite number of powers of prime factors of $R ($by convention a
unit is a least common multiple of the empty family of powers of
prime factors$)$.

{\rm(6)} two ordered monoids $O_1,O_2$ are said to be isomorphic if
there exists a bijection $f$ of $O_1$ onto $O_2$ such that $f$
preserves multiplication and that both $f$ and $f^{-1}$ are isotone
$($in short, $OM$-isomorphic$)$.
\end{defn}

Evidently if $O_1,O_2$ are $OM$-isomorphic, then they are, a
fortiori, order-isomorphic or monoid-isomorphic.

\begin{exam}
{\rm (1)} $G^\ast$ is $OM$-isomorphic to $L^\ast$, the collection of
all nonzero principal ideals of $R$ ordered by inverse inclusion.

{\rm(2)} $G^\ast$ is $OM$-isomorphic to
$M^{\ast\ast\ast}=\{\downarrow[x]|[x]\in G^\ast\}$ endowed with
$\cdot$ defined by $(\downarrow [x])\cdot(\downarrow [y])=\downarrow
([x]\cdot[y])$.
\end{exam}

\begin{thm}
The conditions $D_1, D_2, D_3$ and $D_4$ are equivalent to $D_5$.
\end{thm}
{\bf Proof.} ($D_1$ implies $D_2$) Assume that $J\in
\vartriangle_a$. Since $a=\vee \downarrow a\cap B$, there exists
$b\in \downarrow a\cap B$ such that $b\notin J$ because otherwise $
\downarrow a\cap B\subseteq J$ would imply $a\in J$, a
contradiction. Hence $J\in \vartriangle_b$. Besides, we have
$b\leqslant a$ and so $a\notin J_b$, whence $J_b\in \vartriangle_a$.
it follows that $D_2$ holds.

($D_2$ implies $D_3$) We have $J\in \vartriangle_b$ for each $b\in
B\setminus J$ and hence $J\subseteq J_b$ because $J_b$ is the
greatest $B$-ideal in $\vartriangle_b$. Thus $J\subseteq \cap_{b\in
B\setminus J}J_b$. On the other hand, let $a\in \cap_{b\in
B\setminus J}J_b$. if $a\notin J$, then by $D_2$ there would exist
$b'\in B$ such that $J\in \vartriangle_{b'}$ and $J_{b'}\in
\vartriangle_a$. Since $b'\in B\setminus J, \cap_{b\in B\setminus
J}J_b\subseteq J_{b'}$. But $a\notin J_{b'}$, a fortiori, $a\notin
\cap_{b\in B\setminus J}J_b$, a contradiction. Thus $a\in J$, whence
$\cap_{b\in B\setminus J}J_b\subseteq J$. It follows that
$J=\cap_{b\in B\setminus J}J_b$ and so $D_3$ holds.

($D_3$ implies $D_4$) For each $b\in J\cap B, \downarrow b\subseteq
J$ is clear. Assume that $\downarrow b\subseteq J'\in M$ for all
$b\in J\cap B$. If $J\nsubseteq J'$, then there would exist $a$ such
that $a\in J$ and $a\notin J'$.But by $D_3$ $J'=\cap_{b\in
B\setminus J'}J_b$, so that $a\notin J_{b'}$ for some $b'\in
B\setminus J'$, whence $b'\leqslant a$. Hence $b'\in J$ because $a
\in J$. We have $b'\in J\cap B$ and so $\downarrow b'\subseteq J'$,
whence $b'\in J'$, a contradiction. Thus $J\subseteq J'$. It follows
that $J=\vee_{b\in J\cap b}\downarrow b$ and hence $D_4$ holds.

($D_4$ implies $D_1$) by $D_4$ we have $\downarrow a=\vee_{b\in
\downarrow a\cap B} \downarrow b$. By Example 2.1(2) $G^\ast$ is
$OM$-isomorphic, a fortiori, order-isomorphic to $M^{\ast\ast\ast}$,
whence $a=\vee \downarrow a\cap B$, and so $D_4$ holds.

($D_1$ is equivalent to $D_5$) is trivial.\ \ $\square$

\begin{prop}
{\rm (1)} Under $D_2$ the converse of Proposition 2.1(1) is true.

{\rm (2)} Under $D_3$ the converse of Proposition 2.4 is true.

{\rm (3)} Under $D_3$, $\downarrow b$ with $b\in B$ is strongly
completely $\vee$-irreducible. Its converse is true under $D_4$.
\end{prop}
{\bf Proof.} (1) Assume that $J\in M$. Evidently we have $\uparrow
(G^\ast\setminus J)=G^\ast\setminus J$. Let $a\in G^\ast\setminus
J$, then $a\notin J$ and by $D_2$ there exists $b\in B$ such that
$b\notin J$ and $a\notin J_b$, whence $b\in (G^\ast\setminus J)\cap
B$ and $b\leqslant a$. Hence $G^\ast\setminus J$ is a $B$-set.
Furthermore, suppose $J$ is prime and $a,c\in G^\ast\setminus J$
with $a\wedge c$ existing. If $a\wedge c\notin G^\ast\setminus J$,
then $a\wedge c\in J$ and hence $a\in J$ or $c\in J$, a
contradiction. Thus $a\wedge c\in G^\ast\setminus J$, whence
$G^\ast\setminus J$ is a $B$-filter.

(2)Suppose $J\in M$ is strongly completely $\wedge$-irreducible. By
$D_3$ we have $J=\cap_{b\in B\setminus J}J_b$, whence $J=J_b$ for
some $b\in B\setminus J$.

(3) Assume that $\downarrow b\subseteq \vee J_t$ for any $b\in B$,
$\{J_t\mid t\in T\}\subseteq M$. If $\downarrow b\nsubseteq J_t$ for
each $t$, then $b\notin J_t$ and by $D_3$ we have $J_t=\cap_{b'\in
B\setminus J_t}J_{b'}\subseteq J_b$, whence $\vee J_t\subseteq J_b$,
a contradiction to $b\in \vee J_t$. Thus $\downarrow b\subseteq J_t$
for some $t$ and hence $\downarrow b$ is strongly completely
$\vee$-irreducible. Conversely, suppose $J$ is strongly completely
$\vee$-irreducible. By $D_4$ $J=\vee_{b\in J\cap B} \downarrow b$
and so $J=\downarrow b$ for some $b\in J\cap B$. \ \ $\square$

\begin{rem}\rm
$D_4$ together with Proposition 2.5(3) is equivalent to the fact
that $M$ is order-isomorphic to a complete ring of sets, a fortiori,
$M$ is a completely distributive complete lattice $\cite{Deng}$.
\end{rem}

Now we will turn to uniqueness of arbitrary decomposition. let $P$
be the set of all prime elements of $R$ and $[P]=\{[p]\mid p\in
P\}$, and put $B_{a,[p]}=\{[p]^n\mid B\ni [p]^n\in \downarrow a\}$
for any $[p]\in [P], a\in G^\ast$.

\begin{rem}\rm
If we extend the relation of $a|b$ to $K^\ast=K\setminus \{0\}$,
where $K$ is the quotient field of $R$, then $K^\ast/U$, $U$ being
the set of all units in $R$, is an ordered group and $\leqslant$ is
generated by the integral part $G^\ast$ of $K^\ast/U$. Hence a
result in $K^\ast/U$ also applies to $G^\ast$ if the elements
relating to it belong to $G^\ast$ or it can be expressed in terms of
elements of $G^\ast$. For example, $ab^{-1}\in G^\ast$ with $a,b\in
G^\ast$ can be rephrased like this, $a=bc$ for some $c\in G^\ast$.
On the other hand, a result in a lattice group holds still in
$G^\ast$ if that lattice operation which is carried out on elements
in its condition exists indeed in $G^\ast$.Say, distributive law of
$\cdot$ respect to $\wedge$, i.e., $a\cdot (b\wedge c)=(a\cdot
b)\wedge (a\cdot c)$ can be said to be as follows in $G^\ast$. If
$b\wedge c$ exists, then $(a\cdot b)\wedge (a\cdot c)$ also exists
and is equal to $a\cdot (b\wedge c)$.
\end{rem}

Bearing this remark on mind we will cite some results from
\cite{Bourbaki1} in the remaining of this section and the next when
they are needed and whilst  we will always use notation of
multiplication.

\begin{defn}
{\rm (1)} By the condition $B_1$ we mean that $B_{a,[p]}$ is a
finite set for all $[p]\in [P], a\in G^\ast$.

{\rm(2)} By the condition $B_2$ we mean that $2^{(B)}=2^B$.
\end{defn}

We have $a=\vee \downarrow a\cap B$ under $D_1$. if $B_1$ holds, put
$v_{[p]}(a)=max\{n\mid [p]^n\in \downarrow a\}$, if $B_{a,[p]}\neq
\varnothing$ and $v_{[p]}(a)=0$ if $B_{a,[p]}=\varnothing$ and hence
$a=\vee_{[p]\in [P]}[p]^{v_{[p]}(a)}$.

\begin{prop}
If $D_1,B_1$ hold, then $a=\vee_{[p]\in [P]}[p]^{v_{[p]}(a)}$
uniquely.
\end{prop}
{\bf Proof.} Since $[p]$ is an atom in $G^\ast$ by
(\cite{Bourbaki1}, $\S$1, n$^\circ$13, Proposition 14), $[p]^n\wedge
[q]^m=[1]$ for $[p],[q]\in [P]$ with $[p]\neq [q]$ by
(\cite{Bourbaki1}, $\S$1, n$^\circ$12, Corollary 3$(DIV)$ to
proposition 11$(DIV)$). besides $a=\vee_{[p]\in
[P]}[p]^{v_{[p]}(a)}$ is equivalent to $\downarrow a=\vee \downarrow
[p]^{v_{[p]}(a)}$ by Example 2.1(2). By Proposition
2.5(3)$\downarrow [p]^{v_{[p]}(a)}$ is strongly completely
$\vee$-irreducible, so is $[p]^{v_{[p]}(a)}$ for $A\in 2^{(B)}$.
Thus follows the uniqueness. \ \ $\square$

\begin{prop}
Under $D_1,B_2$ holds if and only if $M$ consists of principal
$B$-ideals only.
\end{prop}
{\bf Proof.} Necessity. Suppose $J\in M$. Then by $B_2$, put $a=\vee
B\cap J$ and we have $B\cap J\subseteq B\cap \downarrow a$. Let
$b\in B\cap \downarrow a$. Since $b\leqslant a=\vee B\cap J$,
$b\leqslant b'$ for some $b'\in B\cap J$ by the last part of the
proof of Proposition 2.6 and hence $b\in J$ because $J$ is a lower
set, whence $b\in B\cap J$. Thus $B\cap \downarrow a\subseteq B\cap
J$ and so $B\cap J=B\cap \downarrow a$. By $D_1$, $J=\vee_{b\in
J\cap B}\downarrow b=\vee_{b\in B\cap \downarrow a}\downarrow
b=\downarrow a$.

Sufficiency. Assume that $M$ consists of principal $B$-ideals only,
then $G$ is order-isomorphic to $M$ by Example 2.1(2), letting $[x]$
correspond to $\downarrow [x]$, in particular, $[0]$ to $\downarrow
[0]$. By Proposition 2.2 $M$ is a complete lattice, so is $G$,
whence $B_2$ holds.\ \ $\square$

We will need the topological representation theorem of \cite{Deng1},
which is concerned in the notion of generalized-continuity of a
poset \cite{Novak}. For convenience we give a direct proof for its
sufficiency.

\begin{lemma}
{\rm \cite{Deng1}} If every element of a complete lattice $P$ can be
decomposed into an arbitrary join of strongly $\vee$-irreducible
elements of $P$, then $P$ is order-isomorphic to the closed-set
lattice ordered by inclusion of some $T_0$-topological space.
\end{lemma}
{\bf Proof.} Let $X$ be the set of all strongly $\vee$-irreducible
elements of $P$ and define $f:P\rightarrow 2^X$ by $f(a)=\downarrow
a\cap X$ for each $a\in P$. It will be readily verified that the
following conditions hold.

(1) $f(0)=\varnothing$ and $f(1)=X$, where $1$ and $0$ denoted the
greatest element and the least element of $P$ respectively.

(2) $f(\vee S)=\cup f(S)$ for any finite subset $S$ of $P$.

(3) $f(\wedge S)=\cap f(S)$ for any $S\subseteq P$.

Let $C=f(P)$ and $C$ satisfies the axioms for closed sets, whence
$X$ endowed with $C$ becomes a topological space.

Consider the corestriction $f^\circ$ of $f$ to $C$ and hence
$f^\circ$ is onto. Since $a=\vee \downarrow a\cap X$, $f^\circ$ is
also one-one and what is more, $a\leqslant b$ is equivalent to
$f^\circ(a)\subseteq f^\circ(b)$. Thus $f^\circ$ is an
order-isomorphism. Besides $(X,C)$ is $T_0$ because
$\{x\}^-=f^\circ(x)$ for $x\in X$, where $\{x\}^-$ is the closure of
$\{x\}$ with $x\in X$, and $f^\circ$ is one-one.\ \ $\square$

\begin{exam}
let $R$ be a commutative ring with identity and $Rd(R)$ the
collection of the radicals of all ideals of $R$. $Rd(R)$ ordered by
inverse inclusion satisfies the condition of Lemma 2.1 and hence has
spec($R$) with Zariski cotopology as its topological representation.
\end{exam}

\begin{prop}
If $D_1$ holds, then $M$ has a topological representation, in which
$\downarrow b$ is the point-closure for each $b\in B$ and all
$J_b$'s with $b\in B$ constitute a base for closed sets.
\end{prop}
{\bf Proof.} By Theorem 2.1 $D_4$ holds and so $J=\vee_{b\in B\cap
J}\downarrow b$, Besides, $\downarrow b$ is strongly completely
$\vee$-irreducible by Proposition 2.5(3),a fortiori, strongly
$\vee$-irreducible, whence $M$ has a topological representation by
Lemma 2.1, of which the proof shows $\downarrow b$ is the
point-closure for any $b\in B$. By $D_3$ $J_b$'s constitute a base
for closed sets.\ \ $\square$

\begin{rem}\rm
By Theorem 2.1 and Propositions 2.1(1), 2.5(1) the collection of
$B$-sets ordered by inclusion is order-isomorphic to the open -set
lattice ordered by inclusion of some $T_0$-topological space, in
which $B$-filter of form $J^c_b$, where $J^c_b$ is the complement of
$J_b$, constitute a base for open sets.
\end{rem}

\section{The fundamental Galois connection}

We will need some knowledge of a Galois connection which is phrased
in the following. Assume that $P_1,P_2$ are posets and $d:
P_1\rightarrow P_2$, $g:P_2\rightarrow P_1$ isotone mappings.

\begin{defn}
{\rm \cite{Novak}} $(g,d)$ is called a Galois connection between
$P_2$ and $P_1$ if $a\leqslant g(b)$ is equivalent to $d(a)\leqslant
b$ for all $a\in P_1, b\in P_2$. And $g$ is called the upper adjoint
of $d$ and $d$ the lower adjoint of $g$.
\end{defn}

\begin{lemma}
{\rm \cite{Novak}(1)} $(g,d)$ is a Galois connection if and only if
$d\circ g(b)\leqslant b$ and $a\leqslant g\circ d(a)$ for all $a\in
P_1, b\in P_2$.

Assume that $(g,d)$ is a Galois connection.

{\rm (2)} $d$ is onto if and only if $g$ is one-one, which in turn
is equivalent to $d\circ g(b)=b$ for any $b\in P_2$.

{\rm (3)} $g$ is onto if and only if $d$ is one-one, which in turn
is equivalent to $g\circ d(a)=a$ for any $a\in P_1$.

{\rm (4)} $d$ preserves existing arbitrary joins and $g$ existing
arbitrary meets.

{\rm (5)} $g(P_2)$ and $d(P_1)$ are order-isomorphic.

{\rm (6)} $d(a)=$min $g^{-1}(\uparrow a)$ and $g(b)=$max
$d^{-1}(\downarrow b)$ for all $a\in P_1, b\in P_2$.
\end{lemma}

\begin{exam}
Assume that $R_1,R_2$ are commutative rings with identity and
$f:R_1\rightarrow R_2$ a homomorphism preserving identity. Let
$Id(R_1), Id(R_2)$ be the collection of all ideals, ordered by
inverse inclusion, of $R_1,R_2$ respectively and define $e:
Id(R_1)\rightarrow Id(R_2)$ by $e(a)=$ the ideal generated by $f(a)$
in $R_2$, for any $a\in Id(R_1)$. Then it will be readily verified
that $(e,f^{-1})$ is a Galois connection between $Id(R_1)$ and
$Id(R_2)$.
\end{exam}

In the following put $I=$ the collection of all ideals of $R$ and
denote its elements by $a,b$. $L$ will be the collection of all
principal ideals. Its elements will be denoted by $(x),(y)$ with
$x,y \in R$, while elements of $G$ by $[x], [y]$ with $x,y \in R$.
by Example 2.1(1) $G$ is $OM$-isomorphic to $L$ ordered by inverse
inclusion and $G^\ast$ to $L^\ast=L\setminus \{(0)\}$.

As was done in \cite{Bourbaki}, for $a\in I^\ast=I\setminus
\{(0)\}$, $a^-$ will denote the divisorial ideal associated with
$a$, i.e., $a^-=\cap \{(x)\mid (x)\supseteq a\}$, and $div(a)$ the
divisor of $a$, i.e., the equivalence class of $a$ (the equivalence
relation is generated by the preorder $\prec$, which is defined as
$a\prec b$ if $\{(x)\mid (x)\supseteq a\}\subseteq \{(x)\mid
(x)\supseteq b\}$). Note that in defining $a^-$ or $\prec$ we use
principal integral ideals only because $a$ is an integral ideal and
so we dispense with principal fractional ideals. We use $D^+$ to
denote the set of all integral divisors. Throughout this treatise we
will adopt notation of multiplication while treating the monoid
structure on divisors.

Put $L^\ast_a=\{(x)\mid L\ni (x)\supseteq a\}$ for any $a\in I$, and
define $g$ by $g(a)=\{[x]\mid (x)\in L^\ast_a\}$ for any $a\in I$
and $d:M\rightarrow I$ by $d(J)=\cap_{[x]\in J}(x)$ for any $J\in
M$. We have $g((0))=\downarrow [0]=G^\ast\in M$ and $d(\downarrow
[0])=(0)$. And what is more, we claim that for any $J$ if
$d(J)=(0)$, then $J=\downarrow [0]$. In fact, $d(J)=\cap_{[x]\in
J}(x)=(0)$ is equivalent to $\vee_{[x]\in J}[x]=[0]$, which in turn
is equivalent to $J=\vee_{[x]\in J}\downarrow [x]=\downarrow [0]$.
Thus put $M^\ast=M\setminus \{\downarrow [0]\}$ and $d$ can be
regarded as a mapping $M^\ast\rightarrow I^\ast$.

Following lemma exhibits $g$ can be also regarded as a mapping $
I^\ast\rightarrow M^\ast$.

\begin{lemma}
$g(a)\in M^\ast$ for any $a\in I^\ast$.
\end{lemma}
{\bf Proof.} Let $[y]\leqslant [x]\in g(a)$, then $(y)\supseteq
(x)\in L^\ast_a$ and so $(y)\in L^\ast_a$, whence $[y]\in g(a)$.
Assume that $2^{(G^\ast)}\ni \downarrow [x]\cap B\subseteq g(a)$.
Let $[z]=\vee \downarrow [x]\cap B$. Then we have $(z)=\cap \uparrow
(x)\cap B^\ast$, where $B^\ast=\{(p^n)\mid [p^n]\in B\}$, and
$\uparrow (x)\cap B^\ast\subseteq L^\ast_a$. Thus $(z)\in L^\ast_a$,
whence $[z]\in g(a)$.\ \ $\square$

\begin{lemma}
$d(J)$ is a divisorial ideal for any $J\in M^\ast$.
\end{lemma}
{\bf Proof.} let $a=d(J)=\cap_{[x]\in J}(x)\neq (0)$. We have
$a\subseteq a^-=\cap \{(x)\mid (x)\supseteq a\}\subseteq
\cap_{[x]\in J}(x)=a$, whence $a^-=a$. Thus $d(J)$ is a divisorial
ideal. \ \ $\square$

~

Let $I^\sim$ be the dual of $I$, i.e., $I$ ordered by inverse
inclusion.

\begin{prop}
$(g,d)$ is a Galois connection between $I^\sim$ and $M$.
\end{prop}
{\bf Proof.} Denote by $\leqslant$ the order of $I^\sim$ and then it
will be readily verified that $d\circ g(a)\leqslant a$ and $g\circ
d(J)\supseteq J$ for any $a\in I^\sim, J\in M$. By Lemma 3.1(1),
$(g,d)$ is a Galois connection.

\begin{defn}
The Galois connection $(g,d)$ between $I^\sim$ and $M$ is called the
fundamental Galois connection.
\end{defn}

\begin{rem}
\rm $(g,d)$ being regarded as Galois connection between
$I^{\ast\sim}$ and $M^\ast$, $d(M^\ast)$ is the collection of all
divisorial ideals of $I^\ast$ (denote it by $D(I^\ast)^\sim$) and is
order-isomorphic to $g(I^{\ast \sim})$ (denote it by $M^{\ast\ast}$)
by Lemma 3.1(5). Thus $D^+$ is also order-isomorphic to
$M^{\ast\ast}$. We can transport the monoid structure on $D^+$ to
$M^{\ast\ast}$ by this isomorphism and then $M^{\ast\ast}$ is an
ordered monoid, whence $D^+$ and $M^{\ast\ast}$ are $OM$-isomorphic.
We choose $D^+$ instead of $D(I^\ast)^\sim$ because $D(I^\ast)^\sim$
is not closed under multiplication. Besides, let $i$ be the
inclusion mapping of $M^{\ast\ast}$ into $M^\ast$. It will be
readily verified that $(i,g\circ d)$ is a Galois connection between
$M^{\ast\ast}$ and $M^\ast$.
\end{rem}

\begin{lemma}
{\rm (1)} $J=\cap_{J\subseteq \downarrow [x]} \downarrow [x]$ for
any $J\in M^{\ast\ast}$.

{\rm (2)} For each $J\in M^\ast$ there exists $J'\in M^{\ast\ast}$
such that $J\subseteq J'$.
\end{lemma}
{\bf Proof.} (1) $J=g(a^-)$ for some $a^-\in I^\ast$. But
$a^-=\vee_{x\in a^-}(x)$ in $I^\ast$, i.e., $a^-=\wedge_{x\in
a^-}(x)$ in $I^{\ast\sim}$, and hence $J=g(a^-)=\cap_{x\in
a^-}g((x))=\cap_{J\subseteq \downarrow [x]}\downarrow [x]$ because
$g$ is $\wedge$-preserving.

(2) $J\subseteq g\circ d(J)$ and $g\circ d(J)\in M^{\ast\ast}$.\ \
$\square$

~

Now we can identify $M^{\ast\ast}$ with $D^+$ and study it in
detail. Let $D$ be the collection of all divisors of $R$.

\begin{lemma}
$M^{\ast\ast}$ is a lattice monoid and the the distributive law of
$\cdot$ with respect to $\wedge$ holds.
\end{lemma}
{\bf Proof.} By (\cite{Bourbaki}, chapter VII, $\S$1, $n^\circ$1,
Theorem 2 (iii) and $n^\circ$2) $D^+$ satisfies the above mentioned
condition, so does $M^{\ast\ast}$ because $M^{\ast\ast}$ is
$OM$-isomorphic to $D^+$. \ \ $\square$

\begin{lemma}
For each $[p]\in [P]$, $\downarrow [p]$ is both a prime element and
an atom of $M^{\ast\ast}$.
\end{lemma}
{\bf Proof.} We need only to show that $div(d(\downarrow
[p]))=div((p))$ is prime in $D^+$, because $M^{\ast\ast}$ is
$OM$-isomorphic to $D^+$. Assume that $div((p))\leqslant
(div(a))\cdot (div(b))=div(a\cdot b)$ with $a,b\in I^\ast$. Then
$(p)\prec a\cdot b$, i.e., $L^\ast_{(p)}\subseteq L^\ast_{a\cdot
b}$. Since $(p)\in L^\ast_{(p)}$, we have $(p)\in L^\ast_{a\cdot b}$
and so $a\subseteq (p)$ or $b\subseteq (p)$, whence $(p)\prec a$ or
$(p)\prec b$. We have $div((p))\leqslant div(a)$ or
$div((p))\leqslant div(b)$ and hence $div((p))$ is prime. Besides,
let $J\in M^{\ast\ast}$ and $J\subseteq \downarrow [p]$. Then there
exists unique $a^-\in D(I^\ast)^\sim$ such that $g(a^-)=J$ so that
$g(a^-)\subseteq g((p))$, whence $a^-\leqslant (p)$ in
$I^{\ast\sim}$, i.e., $a^-\prec (p)$. Now assume that $a^-\subseteq
(x)$ for any $(x)\in I^\ast$, then $(p)\subseteq (x)$ and so
$[x]\leqslant [p]$ in $G^\ast$, whence $[x]=[1]$ or $[x]=[p]$
because $[p]$ is an atom in $G^\ast$ by (\cite{Bourbaki1}, $\S$1,
$n^\circ$13, Proposition 14) in notation of multiplication and
noting that what its proof needs is satisfied although $G^\ast$ is
the integral part of the ordered group, $K^\ast/U$, where $K$ is the
quotient field of $R$, $K^\ast=K\setminus \{0\}$ and $U$ the set of
all units of $R$, only. Thus $g(a^-)\subseteq \{[1],[p]\}$, whence
$J=g(a^-)=\downarrow [1]$ or $J=g(a^-)=\downarrow [p]$. it follows
that $\downarrow p$ is an atom in $M^{\ast\ast}$. \ \ $\square$

\begin{lemma}
If $B_1$ holds, then $B_{J,[p]}=\{\downarrow [p]^n\mid \downarrow
[p]^n\subseteq J,n=1, 2, \cdots\}$ is a finite set for any $[p]\in
[P], J\in M^\ast$.
\end{lemma}
{\bf Proof.} We have $\cap \{(p)^n\mid n=1, 2, \cdots\}=(0)$ because
$B_1$ holds. If $B_{J,[p]}$ were an infinite set, then $J$, being a
lower set, would contain $\{\downarrow [p]^n\mid n=1, 2, \cdots\}$,
whence $\cap_{[x]\in J}(x)\subseteq \cap \{(p)^n\mid n=1, 2,
\cdots\}=(0)$ and hence $\cap_{[x]\in J}(x)=(0)$, a contradiction to
$J\in M^\ast$. Thus $B_{J,[p]}$ is a finite set.\ \ $\square$

~

Assume that $D_1,B_1$ hold. Put $v_{[p]}(J)=max \{n\mid\downarrow
[p]^n\in B_{J,[p]}\}$ if $B_{J,[p]}\neq \varnothing$ and
$v_{[p]}(J)=0$ if $B_{J,[p]}=\varnothing$. $v_{[p]}(J)$ is
well-defined and then $J=\vee_{[p]\in [P]}\downarrow
[p]^{v_{[p]}(J)}$ for $J\in M^\ast$ because $D_1$ is equivalent to
$D_4$ by Theorem 2.1.

\begin{prop}
If $D_1, B_1$ hold, then $J=\vee_{[p]\in [P]}\downarrow
[p]^{v_{[p]}(J)}$ uniquely for any $J\in M^\ast$.
\end{prop}
{\bf Proof.} $\downarrow [p]$ is an atom in $M^{\ast\ast}$ by Lemma
3.6 and so we have $\downarrow [p]^n\wedge \downarrow
[q]^m=\downarrow [1]$ in $M^{\ast\ast}$ for $[p]\neq [q]$ with $[p],
[q]\in [P]$ by (\cite{Bourbaki1}, $\S$1, $n^\circ$ 12, Corollary
$(DIV)$ to Proposition 11 $(DIV)$) because distributive law of
$\cdot$ with respect to $\wedge$ holds in $M^{\ast\ast}$ by Lemma
3.5. Since $(i,g\circ d)$ is a Galois connection between
$M^{\ast\ast}$ and $M^\ast$ by Remark 3.1,$i$ preserves arbitrary
meets, whence $M^{\ast\ast}$ is closed under meets. Thus $\downarrow
[p]^n\wedge \downarrow [q]^m=\downarrow [1]$ also holds in $M^\ast$.
Furthermore, $\downarrow [p]^n$ is strongly completely
$\vee$-irreducible by Proposition 2.5(3). It follows that
$J=\vee_{[p]\in [P]}\downarrow [p]^{v_{[p]}(J)}$ uniquely.\ \
$\square$

\begin{cor}
If the conditions of Proposition 3.2 hold, then each integral
divisor can be decomposed into an arbitrary join of powers of atoms
in $D^+$.
\end{cor}
{\bf Proof.} Note first that $\downarrow [p]$ is an atom in
$M^{\ast\ast}$. Since $D^+$ is $OM$-isomorphic to $M^{\ast\ast}$ and
$B_{J,p}\subseteq M^{\ast\ast}\subseteq M^\ast$, the result follows
from Proposition 3.2. \ \ $\square$

~

In the following we will treat poset-theoretic constructions which
relate decomposition to order representation and structure problems.

Throughout the following we will assume that $D_1, B_1$ hold.

\begin{defn}
Let $P$ be a poset, $P'$ its subposet and $i: P'\rightarrow P$ the
inclusion mapping.

{\rm (1)} $P'$ is said to be of the first kind if $i$ is the lower
adjoint of some Galois connection between $P$ and $P'$.

{\rm (2)} $P'$ is said to be of the second kind if $i$ is the upper
adjoint of some Galois connection between $P'$ and $P$.
\end{defn}

\begin{exam}
{\rm (1)} Put $B_{[p]}^\circ=\{\downarrow [p]^n\mid n=0, 1, 2,
\cdots\}$ for any $[p]\in [P]$. Then $B_{[p]}^\circ$ is a subposet
of $M^\ast$. Put $i_{[q]}: B_{[q]}^\circ\rightarrow M^\ast$ defined
by $i_{[q]}(\downarrow [q]^n)=J$, where $v_{[q]}(J)=n$ and
$v_{[p]}(J)=0$ for $[p]\neq [q]$, and $r_{[q]}: M^\ast\rightarrow
B_{[q]}^\circ$ defined by $r_{[q]}(J)=\downarrow [q]^{v_{[q]}(J)}$.
By Proposition 3.2 it will be readily verified that
$(r_{[q]},i_{[q]})$ is a Galois connection between $M^\ast$ and
$B_{[q]}^\circ$, whence $B_{[p]}^\circ$ is of the first kind.

{\rm (2)} $M^{\ast\ast}$ is a subposet of $M^\ast$ and by Remark 3.1
$i$ is the upper adjoint of Galois connection $(i,g\circ d)$ between
$M^{\ast\ast}$ and $M^\ast$, whence $M^{\ast\ast}$ is of the second
kind. Note that (2) dispenses with $D_1, B_1$.
\end{exam}

\begin{defn}
Let $P$ be a poset with the least element $1$, and $\{P_t\mid t\in
T\}$ a family of subposets of the first kind of $P$ such that each
contains $1$ and that $\vee a_t$ exists uniquely for any $\{a_t\mid
t\in P_t\}$, and assume that $P=\{\vee a_t\mid a_t\in P_t\}$. Then
$P$ is called the internal product of the family. Denote it by
$P=\Pi^iP_t$.
\end{defn}

Let $i_{t'}$ be the inclusion mapping of $P_{t'}$. Evidently we have
$i_{t'}(a_{t'})=\vee b_{t}$, where $ b_{t'}=a_{t'}$ and $b_t=1$ for
$t\neq t'$.

Since $i_{t'}$ is a lower adjoint, its upper adjoint $r_{t'}$
satisfies that $r_{t'}(\vee a_t)=max$ $i_{t'}^{-1}(\downarrow (\vee
a_t))=a_{t'}$ by lemma 3.1 (6), whence $(r_{t'},i_{t'})$ is a Galois
connection between $\Pi^iP_t$ and $P_{t'}$. It is trivial that the
following result holds.

\begin{lemma}
$r_{t'}$ is onto and $a=\vee_{t\in T}i_t\circ r_t(a)$ for any $t'\in
T, a\in \Pi^iP_t$.
\end{lemma}

\begin{exam}
By Example 3.2(1) $B_{[p]}^\circ$ is a subposet of the first kind of
$M^\ast$. By Proposition 3.2 for every $J\in M^\ast$ we have
$J=\vee_{[p]\in [P]}\downarrow [p]^{v_{[p]}(J)}$ uniquely. Thus
$M^\ast=\Pi^i_{[p]\in [P]}B_{[p]}^\circ$. And what is more,
$B_{[p]}^\circ\cap B_{[q]}^\circ=\downarrow [1]$ for $[p]\neq [q]$.
\end{exam}

\begin{defn}
Let $\{P_t\mid t\in T\}$ be a family  of posets containing the least
element $1_t$. The Cartesian product $\Pi P_t$ ordered by the
product order is called the external product of the family. Denote
it by $\Pi^eP_t$.
\end{defn}

Define $r^\ast_{t'}:\Pi^eP_t\rightarrow P_{t'}$ by
$r^\ast_{t'}((a_t))=a_{t'}$ and $i^\ast_{t'}:P_{t'}\rightarrow
\Pi^eP_t$ by $i^\ast_{t'}(a_{t'})=(b_t)$, where $b_{t'}=a_{t'}$ and
$b_t=1_t$ for $t\neq t'$.

\begin{lemma}
$(r^\ast_t, i^\ast_t)$ is a Galois connection, $r^\ast_t$ is onto
and $a=\vee i^\ast_t\circ r^\ast_t(a)$ for any $a\in \Pi^eP_t$.
\end{lemma}
{\bf Proof.} It will readily verified that $r^\ast_t\circ
i^\ast_t(a_t)=a_t$ and $i^\ast_t\circ r^\ast_t(a)\leqslant a$,
whence $(r^\ast_t, i^\ast_t)$ is a Galois connection, $r^\ast_t$ is
onto. It is trivial that $a=\vee i^\ast_t\circ r^\ast_t(a)$.\ \
$\square$

\begin{prop}
$\Pi^iP_t$ is order-isomorphic to $\Pi^eP_t$.
\end{prop}
{\bf Proof.} Define $f:\Pi^iP_t\rightarrow \Pi^eP_t$ by $f(\vee
a_t)=(r_t(\vee a_t))$ for each $\vee a_t\in \Pi^iP_t$. Assume that
$(a_t)\in \Pi^eP_t$, then $a=\vee a_t\in \Pi^iP_t$ with $r_t(a)=a_t$
and so $f(a)=(r_t(a))=(a_t)$, whence $f$ is onto. Furthermore,
$a\leqslant b$ is equivalent to $r_t(a)\leqslant r_t(b)$ for each
$t\in T$, which in turn is equivalent to $(r_t(a))\leqslant
(r_t(b))$, i.e., $f(a)\leqslant f(b)$. It follows that $f$ is an
order-isomorphism.\ \ $\square$

\begin{exam}
$M^\ast$ is order-isomorphic to $\Pi^eB^\circ_{[p]}$, which follows
from Example 3.3 and Proposition 3.3.
\end{exam}

\begin{exam}
$(${\bf Order representation}$)$ Let $N_{[p]}= N=\{0, 1, 2,
\cdots\}$ for each $[p]\in [P]$. Then $N_{[p]}$ ordered by the order
defined by addition is a totally ordered monoid. Put
$f:\Pi^eB^\circ_{[p]}\rightarrow \Pi^eN_{[p]}$ by $f((\downarrow
[p]^{n_{[p]}}))=(n_{[p]})$, then $f$ is a bijection and both $f$ and
$f^{-1}$ are isotone, whence two ordered monoids
$\Pi^eB^\circ_{[p]}$ and $\Pi^eN_{[p]}$  are order-isomorphic. From
Example 3.4 it follows that $M^\ast$ is order-isomorphic to
$\Pi^eN_{[p]}$, which is the order representation of $M^\ast$ under
$D_1, B_1$.
\end{exam}

\begin{exam}
$(${\bf Algebraic interpretation}$)$ As is easily known,
$B^\circ_{[p]}$ is a totally ordered monoid and the order of
$B^\circ_{[p]}$ is defined by the multiplication. Now
$\Pi^eB^\circ_{[p]}$ endowed with componentwise multiplication is an
ordered monoid, of which the product order is also defined by this
multiplication. Through the order-isomorphism of $M^\ast$ onto
$\Pi^eB^\circ_{[p]}$, which follows from Example3.4, we can
transport the monoid structure on $\Pi^eB^\circ_{[p]}$ to $M^\ast$
and hence $M^\ast$ becomes an ordered monoid and $\subseteq$ in
$M^\ast$ is defined by this multiplication. Thus $M^\ast$ is
$OM$-isomorphic to $\Pi^eB^\circ_{[p]}$ and the multiplication
induced on $M^{\ast\ast}$ by this multiplication overlaps with the
original one transported from $D^+$ by $g$ of the fundamental Galois
connection $(g,d)$.

$M^\ast$ as a monoid is monoid-isomorphic to $\Pi^eB^\circ_{[p]}$,
which is the product of a family of monoids in algebraic sense.
\end{exam}

\begin{exam}
$(${\bf Topological interpretation}$)$ By Proposition 2.8 $M$ has a
$T_0$-topological representation. From the proof of Lemma 2.1 $B$
can be taken as the space $X$ and $\{\downarrow b\mid b\in B\}$ is
the collection of all point-closures. Now we consider $\triangle_b$
and put $\triangle_b^c=\{J^c\mid J\in \triangle_b\}$, where we use
$J^c$ to denote the complement of $J$, whence $J^c$ is an open set
in $X$.

Since $J\in \triangle_b$ is equivalent to $b\in J^c, \triangle_b^c$
is the open neighborhood base of $b$ and $J_b^c$ is the least of all
$J^c$'s of $\triangle_b^c$. Furthermore, $D_2$ signifies that for
any $J^c\in \triangle_b^c$ there exist $b'$ such that
$J_{b'}^c\subseteq J^c$ and $J_{b'}^c\in \triangle_b^c$. That is
just the most important one of neighborhood axioms, i.e., for any
neighborhood $N$ of a point $x$ there exists neighborhood $U$ of $x$
with $U\subseteq N$ such that $N$ is a neighborhood of each $y\in
U$. Here $J_b$ is an important tool to study decomposition as is
done by neighborhoods in studying topological local properties.
\end{exam}

\begin{rem}\rm
From Example 3.4 we know that if $J=\vee_{[p]\in [P]}\downarrow
[p]^{v_{[p]}(J)}$ uniquely for each $J\in M^\ast$, then $M^\ast$ is
$OM$-isomorphic to $\Pi^eN_{[p]}$. The converse is true. Let
$1_{[p]}\in \Pi^eN_{[p]}$ such that $r_{[q]}^\ast(1_{[p]})=1$ for
$[q]=[p]$ and $r_{[q]}^\ast(1_{[p]})=0$ for $[q]\neq [p]$. Evidently
$1_{[p]}$ is both prime and an atom. For any $a\in \Pi^eN_{[p]}$, we
have $i_{[p]}^\ast\circ r_{[p]}^\ast(a)=r_{[p]}^\ast(a)\cdot
1_{[p]}$. By lemma 3.9 $a=\vee i_{[p]}^\ast\circ
r_{[p]}^\ast(a)=\vee r_{[p]}^\ast(a)\cdot 1_{[p]}$. Let $f$ be the
$OM$-isomorphism and $J\in M^\ast$ then $J=f(a)$ for unique $a\in
\Pi^eN_{[p]}$ and we have $J=f(a)=\vee
f(1_{[p]})^{r_{[p]}^\ast(a)}$. Uniqueness is clear.
\end{rem}

\begin{center}
College of Mathematics and Econometrics, \\Hunan University,
Changsha 410082, Hunan, China.
\end{center}
\end{document}


\title{\bf Factorization from a poset-theoretic view II
\thanks {Mathematics Subject Classifications(2000): Primary 06A11  Secondary 13A15.}
\thanks {Keywords and phases: $B$-ideals, divisorial
ideals, a Galois connection.}}
\author{\small Zike Deng }
\date{}
\maketitle \noindent{\bf Abstract:} We give several necessary and
sufficient conditions for an element of a monoid to be decomposed
(uniquely) into a product of powers of prime factors in terms of the
properties of $B$-ideals (as well as the properties of powers of
prime factors themselves). As monoids are included a domain and the
collection of integral divisors of a domain.

\section{introduction}
This paper, as the second part of the treatise, continues the study
of \cite{Deng0} and is concerned with finite decomposition. In $\S$2
several necessary and sufficient conditions for finite decomposition
and its uniqueness are obtained (Theorem 2.1, Proposition 2.13).
Internal and external direct products are introduced. $\S$3 gives
several characterizations of a Krull domain (Theorem 3.1). The
mutual relations among principal ideal domains, unique factorization
domains, Dedekind domains and Krull domains are established (Theorem
3.2, Corollary 3.1). Besides, the relation of internal direct
products to subposets of the second kind as well as order
representation and algebraic or topological interpretation relating
to poset-theoretic constructions are obtained in $\S$2 and $\S$3
separately (Remarks 2.1, 3.2, 3.3).

\section{Finite decomposition}
We continue the study of (\cite{Deng0}, $\S$2). $G^\ast, B$ are
defined as before and we use $a, b$ to denote elements of $G^\ast$
for brevity. For other notations the reader is referred to
\cite{Deng0}.

In this section and the next  we will cite some results in
multiplication notation on lattice group from \cite{Bourbaki1} in
case of ordered group (in this section) or lattice monoid (in the
next). Hence those results must be weakened and will be marked with
an asterisk. They need either the weak form of distributive law of
$\cdot$ with respect to $\wedge$ (by $a\cdot (b\wedge c)=(a\cdot
b)\wedge (a\cdot c)$ we mean that if either side exists, so does the
other and both are equal. Denote it by $Dist$) or the fact that
$\leqslant$ is defined by $\cdot$ (i.e., $a\leqslant b$ is
equivalent to $b=a\cdot c$ for some $c$. Denote it by $Defi$). Both
$Dist$ and $Defi$ can be implied by the fact that $\leqslant$ is
generated by the integral part (that part $P$ of a group $G$ such
that $G=P\cdot P^{-1}$).

Now we take some which will be used as examples. In the following
$X$ will be an ordered monoid containing the least element $1$.

\begin{exam}
{\rm (\cite{Bourbaki1}, $\S$1, $n^\circ$12, Proposition 11)}$^\ast$
Assume that $Dist$ holds in $X$.

{\rm (1)} If $x\wedge y=1$, then $x\wedge z=x\wedge(y\cdot z)$ for
any $x, y, z\in X$ (i.e., if either side exists, so does the other
and both are equal).

{\rm (2)} (\rm id. Corollary 1)$^\ast$ If $x\wedge y=1$, and
$x\leqslant y\cdot z$, then $x\leqslant z$.

{\rm (3)} (\rm id. Corollary 3)$^\ast$ If $x_i\wedge y_j=1$ for
$i=1,2,\cdots,n$ and $j=1, 2, \cdots, m$, then $(\pi x_i)\wedge (\pi
y_j)=1$.
\end{exam}
{\bf Proof.} (1) Suppose $x\wedge y=1$, then by $Dist$ we have
$z=z\cdot (x\wedge y)=(z\cdot x)\wedge (z\cdot y)$. If $x\wedge z$
exists, then $x\wedge z$ is a lower bound of $\{x, y\cdot z\}$
because $1\leqslant y, z\leqslant z$ implies $z\leqslant y\cdot z$.
Assume that $w\leqslant x, w\leqslant y\cdot z$. Then we have
$w\leqslant x\cdot z$ from what has just been proved and hence
$w\leqslant z$, whence $w\leqslant x\wedge z$. Thus $x\wedge
z=x\wedge (y\cdot z)$. Similarly for case of $x\wedge (y\cdot z)$
existing.

(2) We have $x\wedge (y\cdot z)=x$ and by (1) $x\wedge z=x$, whence
$x\leqslant z$.

(3) First we claim that $x\wedge y=1, x\wedge z=1$ imply $x\wedge
(y\cdot z)=1$. In fact, by (1), we have $x\wedge z=x\wedge (y\cdot
z)$, whence $x\wedge (y\cdot z)=1$. Then by induction on $m$. \ \
$\square$

\begin{exam}
{\rm (\cite{Bourbaki1}, $\S$1, $n^\circ$12, Proposition 12)}$^\ast$
Assume that $Dist$ holds in $X$.

{\rm (1)} If $x\wedge y=1$ and $x \vee y$ exists, then $x\vee
y=x\cdot y$.

{\rm (2)} If $Defi$ holds in $X$ and $x\wedge y=1$, then $x \vee y$
exists and is equal to $x\cdot y$.
\end{exam}
{\bf Proof.} (1) By $Dist$ we have $x\vee y=(x\vee y)\cdot (x\wedge
y)=[(x\vee y)\cdot x]\wedge [(x\vee y)\cdot y]\geqslant (y\cdot
x)\wedge (x\cdot y)=x\cdot y$. Besides, $x\vee y\leqslant x\cdot y$
as was done in Example 2.1 (1). Hence $x\vee y=x\cdot y$.

{\rm (2)} $x\cdot y$ is an upper bound of $x, y$ as we did in
Example 2.1 (1). Assume that $x\leqslant w, y\leqslant w$. Then
$w=x\cdot z$ for some $z$ by $Defi$. Since $y\leqslant w$ and
$x\wedge y=1$, we have $y\leqslant z$ by Example 2.1 (2). Thus
$x\cdot y\leqslant x\cdot z=w$ and hence $x\vee y=x\cdot y$. \ \
$\square$

\begin{exam}
{\rm (\cite{Bourbaki1}, $\S$1, $n^\circ$13, Proposition 14)}$^\ast$
{\rm(1)} If $Dist$ holds and $X$ is a $\wedge$-semilattice, then
each atom is prime.

{\rm (2)} If cancellation law (i.e., $xz=yz$ implies $x=y$) and
$Defi$ holds, then each prime element is an atom.

{\rm (3)} If cancellation law holds and $x<y$, then $xz<yz$.
\end{exam}
{\bf Proof.} {\rm (1)} Suppose $x$ is an atom and $x\leqslant y\cdot
z$. If $x\nleqslant y$, then $x\wedge y=1$ because $x\wedge y$
exists. By Example 2.(1) we have $x\wedge z=x\wedge(y\cdot z)$.
Since $x\wedge(y\cdot z)$ exists and is equal to $x$, $x\wedge z$
also exists and is equal to $x$, whence $x\leqslant z$. It follows
that $x$ is prime.

{\rm (2)} Assume that $x$ is prime and $y\leqslant x$. Then by
$Defi$ we have $x=y\cdot z$ for some $z$ and hence $x\leqslant y$ or
$x\leqslant z$, whence in the former case $x=y$ and in the latter
case $z=x$ because $z\leqslant x$ as we did in Example 2.1(1), and
so $x=y\cdot x$, which implies $y=1$ by cancellation law. Thus $x$
is an atom.

{\rm (3)} Assume that $x<y$. Then we have $xz\leqslant yz$. If
$xz=yz$, by cancellation law $x=y$, a contradiction. Hence $xz<yz$.
\ \ $\square$

~

Now we turn to the properties of $B$. $G^\ast$ is the integral part
of ordered group $K^\ast/U$ and $\leqslant$ is generated by $G^\ast$
so that $Defi$, $Dist$ and cancellation law in Examples 2.1, 2.2,
2.3 all hold. By Example 2.3 (2), $[p]\in B$ with $p\in P$ is an
atom. A finite set of $B$ as such $\{[q_i]^{n_i}\mid [q_i]\neq
[q_j]$ for $i\neq j$ and $i=1, 2, \cdots, m\}$ is called a condensed
set. Any finite set $A$ of $B$ can be reduced to a condensed set
$A^\ast$ such that $A^\ast\subseteq A$ by combining powers of the
same base (taking maximum of the indexes). It is trivial that if
$\vee A^\ast$ exists, then $\vee A$ exists and is equal to $\vee
A^\ast$.

\begin{lemma}
$B$ satisfies the conditions. {\rm (1)} $\downarrow b\cap B$ is
finite set for any $b\in B$. (Denote it $B_3$).

{\rm (2)} $\vee A$ exists for any finite set of $B$ and if $A$ is
condensed, then $\vee A=\prod A$.

{\rm (3)} $a\leqslant \vee A$ implies $a\leqslant b$ for some $b\in
A$, for any $a\in B$, a finite set $A$ with $\vee A$ existing and
$A\subseteq B$. (Denote it $B_4$).
\end{lemma}
{\bf Proof.} (1) Let $b=[p]^n$. Since $q\in B$ is both prime and an
atom, $[q]^m\leqslant [p]^n$ imolies $q=p$ and $m\leqslant n$,
whence $\downarrow b\cap B$ is a finite set.

(2) Assume first that $A=\{[q_i]^{n_i}\mid [q_i]\neq [q_j]$ for
$i\neq j$ and $i=1, 2, \cdots, m\}$ is a condensed set. By Example
2.1(3) we have $[q_i]^{n_i}\wedge [q_j]^{n_j}=[1]$, whence $\vee
[q_i]^{n_i}=\prod [q_i]^{n_i}$ by Example 2.2(2) and by induction on
$m$ noting Example 2.1(3). Next suppose $A$ is any finite set of
$B$. Then there is a condensed set $A^\ast$ such that
$A^\ast\subseteq A$. We have already proved that $\vee A^\ast= \prod
A$, and hence $\vee A$ exists and is equal to $\vee A^\ast$.

(3) If $A$ is condensed, by (2) $\vee A=\prod A$. Since $[p]\in B$
is both prime and an atom, $a\leqslant \prod A$ implies $a\leqslant
b$ for some $b\in A$. Now suppose $A$ is a finite set of $B$. Then
$\vee A=\vee A^\ast$ for some condensed set $A^\ast$ and so
$a\leqslant b$ for some $b\in A^\ast$, a fortiori, $b\in A$ because
$A^\ast\subseteq A$. \ \ $\square$

Now we start to study finite decomposition. For any $C\subseteq
G^\ast$ we use $C^{(f)}$ to denote the collection $\{A\mid A$ is a
finite subset of $C$ such that $\vee A$ exists $\}$.

\begin{defn}
{\rm (1)} By the condition $F_1$ we mean that for each $a\in G^\ast$
there exists $A\in (\downarrow a\cap B)^{(f)}$ such that $a=\vee A$.

{\rm (2)} By the condition $F_2$ we mean that for any $a, b\in
G^\ast, \downarrow b\cap B\in 2^{(B)}, a\leqslant \vee \downarrow
b\cap B$ implies $a\leqslant \vee A$ for some $A\in (\downarrow
b\cap B)^{(f)}$.

{\rm (3)} By $DCC$ we mean that for any descending chain
$a_1\geqslant a_2\geqslant \cdots\geqslant a_n\geqslant \cdots$ in
$G^\ast$ there exists $m$ such that $a_i=a_m$ for $i\geqslant m+1$.
\end{defn}

\begin{prop}
{\rm (1)} $F_1$ is equivalent to $F_2$ and $D_1$.

{\rm (2)} Under $B_3$ $F_1$ is equivalent to $DCC$ and
$ir(G^\ast)\subseteq B$, where $ir(G^\ast)$ is the set $\{a\mid a\in
G^\ast$ such that $a=\vee A$ with $A\in G^{\ast (f)}$ implies $a=b$
for some $b\in A\}$.
\end{prop}
{\bf Proof.} (1) ($F_2$ and $D_1$ imply $F_1$) By $D_1$ we have
$a=\vee \downarrow a\cap B$. By $F_2$ $a\leqslant \vee A$ for some
$A\in (\downarrow a\cap B)^{(f)}$ but $\vee A\leqslant \vee
\downarrow a\cap B=a$, whence $a=\vee A$.

($F_1$ implies $F_2$ and $D_1$) We have $a=\vee A$ for some $A\in
(\downarrow a\cap B)^{(f)}$ by $F_1$. Since $A\subseteq \downarrow
a\cap B, \vee \downarrow a\cap B$ exists and is equal to $a$, whence
$D_1$ holds. Furthermore, let $a\leqslant \vee \downarrow b\cap B\in
2^{(B)}$, then $\vee \downarrow b\cap B=b=\vee A$ for some $A\in
(\downarrow b\cap B)^{(f)}$ and hence $a\leqslant \vee A$. Thus
$F_2$ holds.

(2) ($F_1$ together with $B_3$ implies $DCC$ and
$irr(G^\ast)\subseteq B$) Assume that $a_1>a_2>\cdots>a_n>\cdots$ is
a strictly descending chain in $G^\ast$. By $F_1$ we have $a_n=\vee
A_n$ for some $A_n\in (\downarrow a_n\cap B)^{(f)}$. Let
$a_1=b_1\vee b_2\vee \cdots\vee b_m$ with $\{b_i\mid i=1, 2, \cdots,
m\}=A_1\in (\downarrow a_1\cap B)^{(f)}$. Since $a_n<a_1$ for $n>1$,
by [1,Proposition 2.5(3) and Example 2.1(2)] for any $b\in A_n$ with
$n>1$ there exists $j$ such that $b\leqslant b_j$ and hence $b\in
\downarrow b_j\cap B\subseteq \cup \{\downarrow b_i\cap B\mid i=1,
2, \cdots, m\}$, which is a finite set by $B_3$. Thus $A_n\subseteq
\cup \downarrow b_i\cap B$ for each $n$. Furthermore we have
$A_n\neq A_m$ for $n\neq m$ so that this chain must terminate at
some $n$. On the other hand, suppose $a\in ir(G^\ast)$, then by
$F_1$ we have $a=\vee A$ for some $A\in (\downarrow a\cap B)^{(f)}$
and so $a=b$ for some $b\in A$, whence $a\in B$. Thus
$ir(G^\ast)\subseteq B$.

($DCC$ and $irr(G^\ast)\subseteq B$ imply $F_1$) Put $X=\{a\in
G^\ast\mid a$ can not be decomposed into a finite join of elements
$\in ir(G^\ast)\}$. If $X\neq \varnothing$, then there exists
minimal $m\in X$. Evidently $m\notin ir(G^\ast)$ and hence $m=\vee
a_i$ for some $\{a_i\mid i=1, 2, \cdots, n\}\in G^{\ast (f)}$ with
$a_i<m$. We have $a_i\notin X$, whence $a_i=\vee b^i_j$ for some
$\{b^i_j\mid j=1, 2, \cdots, n_i\}\in (ir(G^\ast))^{(f)}$. Thus
$\vee_{i,j}b^i_j$ exists and is equal to $m$, a contradiction. It
follows that $X=\varnothing$ and hence each $a\in G^\ast$ can be
decomposed into a finite join of elements of $ir(G^\ast)$, which
together with $ir(G^\ast)\subseteq B$ implies that $F_1$ holds. \ \
$\square$

\begin{exam}
Let $R$ be a Noetherian ring and $Id(R)$ the collection of ideals of
$R$, ordered by inverse inclusion. As is well known, $DCC$ holds in
$Id(R)$ and any $\vee$-irreducible ideal is primary. Hence each
ideal can be decomposed into a finite join of primary ideals by
Proposition 2.1(2).
\end{exam}

\begin{defn}
By the condition $F_3$ we mean that for any $a\in G^\ast, \{J_t\mid
t\in T\}\subseteq M$, $\downarrow a\subseteq \vee J_t$ implies
$\downarrow a\subseteq \vee J_{t_i}$ for some $\{t_i\mid
i=1,2,\cdots, n\}\subseteq T$.
\end{defn}

\begin{prop}
$D_1$ and $F_3$ are equivalent to $F_1$.
\end{prop}
{\bf Proof.} ($D_1$ and $F_3$ imply $F_1$) Assume that $a\in
G^\ast$. By (\cite{Deng0}, Theorem 2.1) $D_1$ implies $D_4$ and
hence $\downarrow a=\vee_{b\in \downarrow x\cap B} \downarrow b$,
whence $\downarrow a\leqslant \vee \downarrow b_i$, where $\{b_i\mid
i=1, 2, \cdots, n\}\subseteq \downarrow a\cap B$, by $F_3$. Since
$\vee \downarrow b_i\leqslant \vee \downarrow b=\downarrow a$,
$\downarrow a=\vee \downarrow b_i$, whence by (\cite{Deng0}, Example
2.1(2)) we have $a=\vee b_i$ and hence $F_1$ holds.

($F_1$ implies $D_1$ and $F_3$) For any $a\in G^\ast$ we hace
$a=\vee b_i$ for some $\{b_i\mid i=1, 2, \cdots, n\}\in (\downarrow
a\cap B)^{(f)}$ by $F_1$, whence $\vee \downarrow a\cap B$ exists
and is equal to $a$. Thus $D_1$ holds. Furthermore, assume that
$\downarrow a\subseteq \vee J_t$ for $\{J_t\mid t\in T\}\subseteq
M$. From above we know that $a=\vee b_i$ and by (\cite{Deng0},
Example 2.8 (2)) $\downarrow a=\vee \downarrow b_i$. Thus we have
$\downarrow b_i\subseteq \vee J_t$ and by (\cite{Deng0}, Theorem 2.1
and Proposition 2.5(3)) $\downarrow b_i\subseteq J_{t_i}$ for some
$t_i$, whence $\downarrow b_i\subseteq \vee \{J_{t_j}\mid j=1, 2,
\cdots, n\}$ for each $i$. Hence $\downarrow a\subseteq \vee
\{J_{t_j}\mid j=1, 2, \cdots, n\}$ and so $F_3$ holds. \ \ $\square$

\begin{prop}
If $F_3$ holds, then for all $a\in G^\ast$, $J\in \Delta_a$ there
exists $K\in \Delta_a$ such that $J\subseteq K$ and $K$ is maximal
in $\Delta_a$.
\end{prop}
{\bf Proof.} Let $J\in \Delta_a$ and put $\Delta=\{J'\mid J'\in
\Delta_a$ and $J\subseteq J'\}$. then $\Delta\neq \varnothing$
because $J\in \Delta$. Suppose $\{J'_t\mid t\in T\}\subseteq \Delta$
is a chain. We claim that $\vee J'_t\in \Delta$. In fact,
$J\subseteq \vee J'_t$ is trivial. If $a\in \vee J'_t$, then
$\downarrow a\subseteq \vee J'_t$ and by $F_3$ $\downarrow
a\subseteq \vee J'_{t_i}=max$ $\{J'_{t_i}\}$ for some $\{t_i\mid
i=1, 2, \cdots, n\}\subseteq T$, whence $a\in max$ $\{J'_{t_i}\}\in
\Delta_a$, a contradiction. Hence $\vee J'_t\in \Delta_a$. It
follows that $\vee J'_t\in \Delta$. By Zorn's Lemma there exists a
maximal $K\in \Delta$ and $K$ is clearly maximal in $\Delta_a$. \ \
$\square$

Denote by $\sum_2$ the collection of all such $K'$s. From above we
know that for each $a, J, a\notin J$ implies $J\subseteq K$ for some
$K\in \sum_2$ and denote such a $K$ by $K_{J,a}$.

\begin{prop}
If $F_3$ holds, then $J=\cap_{a\notin J}K_{J,a}$ for each $J\in M$
and $\sum_2=$ the collection of all completely $\wedge$-irreducible
$B$-ideals.
\end{prop}
{\bf Proof.} Given $J\in M$ and by Proposition 2.3 for each $a\notin
J$ there exists $K_{J,a}$ such that $J\subseteq K_{J,a}$. We claim
that $J=\cap_{a\notin J} K_{J,a}$. In fact, $J$ is a lower bound of
$\{K_{J,a}\mid a\notin J\}$. Suppose $J'\subseteq K_{J,a}$ for all
$a\notin J$. If $J'\nsubseteq J$, there would exist $b\in
J'\setminus J$ and so $b\in J'\subseteq \cap_{a\notin J}
K_{J,a}\subseteq K_{J,b}$, a contradiction. Thus $J'\subseteq J$,
whence $J=\cap_{a\notin J} K_{J,a}$.

On the other hand, suppose $\{J_t\mid t\in T\}\subseteq M$ and
$K_{J, a}=\cap J_t$. Since $a\notin K_{J, a}$, we have $a\notin
J_{t_0}$ for some $t_0\in T$ and hence $J_{t_0}\in \Delta_a$, whence
$K_{J_{t_0}, a}\supseteq J_{t_0}$. But $K_{J, a}\subseteq J_{t_0}$
so that $K_{J, a}\subseteq K_{J_{t_0}, a}$. By maximality of $K_{J,
a}$ in $\Delta_a$ we have $K_{J, a}=K_{J_{t_0}, a}$, whence
$J_{t_0}=K_{J, a}$. Thus $K_{J, a}$ is completely
$\wedge$-irreducible. Conversely let $J$ be completely
$\wedge$-irreducible. we have $J=\cap_{a\notin J} K_{J, a}$ and
hence $J=K_{j, a}$ for some $a\notin J$. It follows that $J\in
\sum_2$. \ \ $\square$

\begin{prop}
Under $F_3$, $D_1$ is equivalent to $\sum_2\subseteq \sum_1$.
\end{prop}
{\bf Proof.} ($D_1$ implies $\sum_2\subseteq \sum_1$) Assume that
$D_1$ holds and $K\in \sum_2$. Then $K\in \Delta_a$ for some $a\in
G^\ast$ and so there exits $J_b$ with $b\in B$ such that $K\subseteq
J_b\in \Delta_a$ by $D_2$ because $D_1$is equivalent to $D_2$ by
(\cite{Deng0}, Theorem 2.1). Thus $K=J_b$ by maximality of $K$ in
$\Delta_a$. It follows that $\sum_2\subseteq \sum_1$.

($\sum_2\subseteq\sum_1$ implies $D_1$) Suppose $a\notin J$. by
Proposition 2.3 we have $J\subseteq K\in \Delta_a$ for $K$ maximal
in $\Delta_a$. Since $K\in \sum_2\subseteq \sum_1$, $K=J_b$ for some
$b\in B$ and hence $J\subseteq J_b\in \Delta_a$, which implies $J\in
\Delta_b$. Thus $D_2$ holds, and so $D_1$ also holds because $D_1$
is equivalent to $D_2$ by (\cite{Deng0}.\, Theorem 2.1). \ \
$\square$

\begin{defn}
{\rm (1)} $a, b\in G^\ast$ with $a\neq b$ are said to be
incomparable if neither $a<b$ nor $b<a$.

{\rm (2)} $A\subseteq G^{\ast(f)}$ is said to be pairwise
incomparable if for all distinct $a, b\in A$, $a,b$ are
incomparable.
\end{defn}

\begin{lemma}
Let $a=\vee A$ with $A\in G^{\ast(f)}$. Then there exists a pairwise
incomparable $C\subseteq A$ such that $\vee A=\vee C$.
\end{lemma}
{\bf Proof.} The verification will be easily completed. \ \
$\square$

\begin{prop}
Assume that $F_3$ holds and $\sum_2\subseteq \sum_1$. Then every
$a\in G^\ast$ can be decomposed into a finite join of pairwise
incomparable elements of $B$.
\end{prop}
{\bf Proof.} Given $a\in G^\ast$. By Proposition 2.5 $D_1$ holds and
by $D_1$ we have $\downarrow a=\vee_{b\in B\cap \downarrow a}
\downarrow b$, whence $\downarrow a\leqslant \vee \downarrow b_i$
for some $\{b_i\mid i=1, 2, \cdots, n\}\subseteq B\cap \downarrow a$
by $F_3$. Since $\vee \downarrow b_i\leqslant \vee \downarrow
b=\downarrow a$, $\downarrow a=\vee \downarrow b_i$. Besides,
$\downarrow b_i\leqslant \downarrow a$ and so $b_i\leqslant a$,
which implies $a\notin J_{b_i}$. We claim that for any $J_b$ with
$b\in B$, $a\notin J_b$ implies $b_i\notin J_b$ for some $i$. In
fact, if $b_i\in J_b$ for each $i$, then $a\in \vee \downarrow
b_i\subseteq J_b$, a contradiction.

We will use the following fact, which is easily verified. For any
$a,b\in B$, $ a\leqslant b$ is equivalent to $\Delta_a\subseteq
\Delta_b$, which in turn is equivalent to $J_a\subseteq J_b$. Now
since $\sum_2\subseteq \sum_1$,for each $K_{J, a}$ we have $K_{J
,a}=J_b$ for some $b\in B$, whence $a\notin J_b$. Then $b_i\notin
J_b$ for some $i$ and hence $b\leqslant b_i$, which is equivalent to
$J_b\subseteq J_{b_i}$. Thus $K_{J ,a}\subseteq J_{b_i}$ and by
maximality of $K_{J, a}$ in $\Delta_a$ we have $K_{J, a}=J_{b_i}$.

Put $L=\{b_i\mid J_{b_i}=K_{J, a}$ for some $K_{J, a}\in
\Delta_a\}$. Evidently $L\subseteq \{b_i\mid i=1, 2, \cdots, n\}$.
Now we show that $a=\vee_{b_i\in L} b_i$. $a$ is clearly an upper
bound of $L$. Suppose $b_i\leqslant c$ for each $b_i\in L$, then
$\Delta_{b_i}\subseteq \Delta_c$. We claim that $K_{J, a}\in
\Delta_c$ for any $K_{J, a}\in \Delta_a$. In fact, since $K_{J,
a}=J_{b_i}$ for some $b_i\in L$, we have $K_{J, a}\in \Delta_{b_i}$
and hence $K_{J, a}\in \Delta_c$. Consequently $\Delta_a\subseteq
\Delta_c$ because for any $J\in \Delta_a$ we have $J\subseteq K_{J,
a}$ and $K_{J, a}\in \Delta_c$, a fortiori, $J\in \Delta_c$, whence
$a\leqslant c$. It follows that $a=\vee_{b_i\in L} b_i$.

On the other hand, if $b_i, b_j\in L$ with $b_i\neq b_j$, then
$J_{b_i}=K_{J, a}, J_{b_j}=K_{J', a}$ are maximal in $\Delta_a$.
Hence they are incomparable, so are $b_i, b_j$. \ \ $\square$

\begin{defn}
 By uniqueness of finite decomposition we mean that $\vee A=\vee C$
implies that $A=C$ for all pairwise incomparable $A, C\in B^{(f)}$.
\end{defn}

\begin{prop}
Uniqueness of finite decomposition is equivalent to $B_4$.
\end{prop}
{\bf Proof.} (Uniqueness of finite decomposition implies $B_4$)
Assume that $a\leqslant \vee A$ for $a\in B$, $A\in B^{(f)}$. Then
by Lemma 2.2 there exists a pairwise incomparable $A^\ast\subseteq
A$ such that $\vee A^\ast=\vee A$. Suppose that $a\nleqslant b$ for
every $b\in A^\ast$. If $a, b$ are incomparable for each $b\in
A^\ast$, then $\{a\}\cup A^\ast$ is pairwise incomparable and $\vee
(\{a\}\cup A^\ast)=\vee A^\ast$ and hence by uniqueness of finite
decomposition $\{a\}\cup A^\ast=A^\ast$, whence $a\in A^\ast$, a
contradiction. If $a>b$ for some $b\in A^\ast$, then delete all
those $b'$s and denote by $C$ the set of remaining elements of
$A^\ast$ so that $\{a\}\cup C$ is pairwise incomparable. We have
$\vee (\{a\}\cup C)=\vee A^\ast$ and hence $\{a\}\cup C=A^\ast$,
which implies $a\in A^\ast$, another contradiction. Therefore
$a\leqslant b$ for some $b\in A^\ast\subseteq A$. It follows that
$B_4$ holds.

($B_4$ implies uniqueness of finite decomposition) Suppose $\vee
A=\vee C$ for pairwise incomparable $A,C\in B^{(f)}$. Let $a\in A$,
then $a\leqslant \vee C$ and hence $a\leqslant b$ for some $b\in C$.
For this $b$ in turn there is $c\in A$ such that $b\leqslant c$,
whence $a\leqslant c$. By incomparability of $a, c$ we have $a=c$.
Hence $a=b\in C$. The converse inclusion can be proved similarly. \
\ $\square$

\begin{lemma}
If $F_1$ holds, then $G^\ast$ is a lattice.
\end{lemma}
{\bf Proof.} Let $\{a_i\mid i=1, 2, \cdots, n\}\subseteq G^\ast$ and
we claim that $\vee a_i$ exists. In fact, for each $a_i$ we have
$a_i=\vee A_i$ for some $A_i\in (\downarrow a_i\cap B)^{(f)}$ by
$F_1$ and $\vee (\cup A_i)$ exists by Lemma 2.1(2) so that $\vee
a_i$ exists and is equal to $\vee (\cup A_i)$. We know that $\wedge
a_i$ also exists from the proof of (\cite{Bourbaki1}, $\S$1,
n$^\circ$9, Proposition 8) in notation of multiplication and noting
that $K^\ast/ U$ is an ordered group. \ \ $\square$

\begin{thm} The following conditions except (4) are equivalent.

{\rm (1)} Each $a\in G^\ast$ can be decomposed uniquely into the
product of powers of atoms in $G^\ast$.

{\rm (2)} $F_1$ holds.

{\rm (3)} $F_2$ and one of $D_i$, $i=1, 2, 3, 4$ hold.

{\rm (4)} $DCC$ and $ir(G^\ast)\subseteq B$ imply $F_1$, and $F_1$
together with $B_3$ implies $DCC$ and $ir(G^\ast)\subseteq B$.

{\rm (5)} $F_3$ and one of $D_i$, $i=1, 2, 3, 4$ hold.

{\rm (6)} $F_3$ holds and $\sum_2\subseteq \sum_1$.

{\rm (7)} $R$ is a unique factorization domain.
\end{thm}
{\bf Proof.} (1) implies (2). This follows from Lemma 2.1 (2).

(2) is equivalent to (3). This follows from Proposition 2.1 (1) and
(\cite{Deng0}, Theorem 2.1).

(4) follows from the proof of the Proposition 2.1 (2).

(2) is equivalent to (5). This follows from Proposition 2.2 and
(\cite{Deng0}, Theorem 2.1).

(5) is equivalent to (6). This follows from Proposition 2.5 and
(\cite{Deng0}, Theorem 2.1).

(6) implies (1). This follows from Propositions 2.6, 2.7 and Lemmas
2.1(2)(3), 2.3.

(7) is equivalent to (1). This is trivial, noting that a unit is the
product of empty family of powers of atoms. \ \ $\square$

\begin{prop}
If $F_1$ holds, then $G^\ast$ is a distributive lattice and $B=$ the
set of all strongly $\vee$-irreducible elements of $G^\ast$.
\end{prop}
{\bf Proof.} By Lemma 2.3, $G^\ast$ is a lattice. Since $a\wedge
(b\vee c)\geqslant (a\wedge b)\vee (a\wedge c)$ is trivial, we need
only to show that $a\wedge (b\vee c)\leqslant (a\wedge b)\vee
(a\wedge c)$. For each $d\in B\cap \downarrow (a\wedge (b\vee c))$
we have $d\leqslant a$, and $d\leqslant b\vee c$, whence by Lemma
2.1 (3) $d\leqslant b$ or $d\leqslant c$. Thus $d\leqslant a\wedge
b$ or $d\leqslant a\wedge c$, so that $d\leqslant a\wedge b)\vee
(a\wedge c)$. Since $F_1$ holds, a fortiori, $D_1$ also holds,
$a\wedge (b\vee c)=\vee B\cap \downarrow (a\wedge (b\vee
c))\leqslant (a\wedge b)\vee (a\wedge c)$.

On the other hand, suppose $b\leqslant a\vee c$ for any $b\in B$,
$a, c\in G^\ast$. By $F_1$ $a=\vee A$ with $A\in (\downarrow a\cap
B)^{(f)}$ and $c=\vee C$ with $C\in (\downarrow c\cap B)^{(f)}$ and
from the proof of Lemma 2.3 we know that $a\vee c=\vee (A\cup C)$.
By Lemma 2.1 (3), $b\leqslant d$ for some $d\in A\cup C$, whence
$b\leqslant a$ or $b\leqslant c$. Thus $b$ is strongly
$\vee$-irreducible. Furthermore, let $a$ be  strongly
$\vee$-irreducible. By $F_1$ we have $a=\vee A$ for some $A\in
(\downarrow a\cap B)^{(f)}$, whence $a\leqslant b$ foe some $b\in
A$. Since $b\leqslant a$ is clear, $a=b\in B$. \ \ $\square$

\begin{prop}
If $F_1$ holds, then $M$ is the collection of all lattice ideals of
$G^\ast$.
\end{prop}
{\bf Proof.} Suppose $J$ is a lattice ideal and $2^{(G^\ast)}\ni
\downarrow a\cap B\subseteq J$. By $F_1$ $a=\vee \downarrow a\cap B=
\vee A$, where $A\in (\downarrow a\cap B)^{(f)}$ because, a
fortiori, $D_1$ also holds. Since $A\subseteq \downarrow a\cap
B\subseteq J$, $a=\vee A=\vee \downarrow a\cap B\in J$, whence $J\in
M$. Conversely assume that $J\in M$ and $\{a_1,a_2\}\subseteq J$, by
$F_1$ we have $a_1=\vee A_1$ with $A_1\in (\downarrow a_1\cap
B)^{(f)}$ and $a_2=\vee A_2$ with $A_2\in (\downarrow a_2\cap
B)^{(f)}$.From the proof of Lemma 2.3 we know that $a_1\vee a_2=\vee
(A_1\cup A_2)=\vee \downarrow (a_1\vee a_2)\cap B$. Since $A_1\cup
A_2\subseteq J$ because $J$ is a lower set, we have $a_1\vee a_2\in
J$. Hence $J$ is a lattice ideal. \ \ $\square$

~

Now we turn to the continuing study of (\cite{Deng0}, $\S$3) Assume
that $F_1$ holds, then by Theorem 2.1 (1) and Lemma 2.1 (3) for any
$a\in G^\ast$ we have $a=\vee_{[p]\in [P]}[p]^{v_[p](a)}$ uniquely
with all $v_{[p]}(a)=0$ except for a finite number of them, where
$v_[p](a)=n$ with $[p]^n\in A$, the pairwise incomparable subset of
$B^{(f)}$ corresponding to $a$ and $v_[p](a)=0$ if $[p]$ is not a
base of any power of $A$, if we use join instead of product.

\begin{defn}
{\rm (1)} The subposet of the internal product $\prod^i P_t$ of
$\{P_t\mid t\in T\}$ such that for each $\vee a_t$ with $a_t\in
P_t$, all $a_t$ are the common least element $1$ except for a finite
number of them is called internal direct product of $\{P_t\mid t\in
T\}$. Denote it by $\prod^{id}P_t$.

{\rm (2)} The subposet of the external product $\prod^e P_t$ of
$\{P_t\mid t\in T\}$ such that for each $(a_t)$ with $a_t\in P_t$,
all $a_t$ are the least elements $1_t'$s except for a finite number
of them is called external direct product of $\{P_t\mid t\in T\}$.
Denote it by $\prod^{ed} P_t$.
\end{defn}

Note that $i_t(a_t)\in \prod^{id}P_t$ and $i^\ast_t(a_t)\in
\prod^{ed}P_t$ for any $a_t\in P_t$.

\begin{prop}
Assume that $F_1$ holds. Then $G^\ast=\prod^{id}Q_{[p]}$, where
$Q_{[p]}=\{[p]^n\mid n=0, 1, 2, \cdots\}$, and $\prod^{id}Q_{[p]}$
is order-isomorphic to $\prod^{ed}Q_{[p]}$.
\end{prop}
{\bf proof.} $G^\ast \subseteq \prod^{id}Q_{[p]}$ by what has just
been said above and $\prod^{id}Q_{[p]}\subseteq G^\ast$ by Lemma 2.1
(2). The proof of the last phrase is similar to that of
(\cite{Deng0}, Proposition 3.3). As for $Q_{[p]}$ being a subposet
of the first kind, refer to (\cite{Deng0}, Example 3.2 (1)). \ \
$\square$

~

Since each $Q_{[p]}$ is an submonoid of $G^\ast$, we can endow
$\prod^{id}Q_{[p]}$ with multiplication defined by $a\cdot b=\vee
i_{[p]}(r_{[p]}(a)\cdot r_{[p]}(b))$ for any $a, b\in
\prod^{id}Q_{[p]}$ and it is trivial that $\prod^{id}Q_{[p]}$ is an
ordered monoid. Then $G^\ast$ as a monoid is the internal direct
product of monoids $Q_{[p]}'$s. The same is also true for
$\prod^{ed}Q_{[p]}$ if we endow  it with componentwise
multiplication. Thus it is easily verified that $\prod^{id}Q_{[p]}$
is $OM$-isomorphic to $\prod^{ed}Q_{[p]}$. Furthermore, $N=\{0, 1,
2, \cdots\}$ regarded as a monoid of addition is a totally ordered
monoid if it is ordered by the order defined by addition. Put
$N_{[p]}=N$ for each $[p]\in [P]$ and then $\prod^{ed}N_{[p]}$
becomes an ordered monoid.

\begin{prop}
If $F_1$ holds, then

{\rm (1)} $G^\ast$ is a lattice.

{\rm (2)} $G^\ast$ is $OM$-isomorphic to $\prod^{ed}N_{[p]}$
\end{prop}
{\bf Proof.} (1) This is another proof of a part of Lemma 2.3. We
still denote the Galois connection between $G^\ast$ and $Q_{[p]}$ by
$(r_{[p]}, i_{[p]})$. Each $Q_{[p]}$ is a totally ordered set and so
is a lattice. Besides, we have $a=\vee i_{[p]}\circ r_{[p]}(a)$. As
is easily verified, for $a, b\in G^\ast$ we have $a\vee b=\vee
i_{[p]}(r_{[p]}(a)\vee r_{[p]}(b))$ and $a\wedge b=\vee
i_{[p]}(r_{[p]}(a)\wedge r_{[p]}(b))$, see (\cite{Deng0}, lemma
3.8).

(2) Because $Q_{[p]}=\{[p]^n\mid n=0, 1, 2, \cdots\}$
is$OM$-isomorphic to $N_{[p]}=\{0, 1, 2, \cdots\}$, we have that
$\prod^{id}Q_{[p]}$ is $OM$-isomorphic to $\prod^{ed}Q_{[p]}$, which
in turn is $OM$-isomorphic to $\prod^{ed}N_{[p]}$. \ \ $\square$

Now we turn to finite decomposition and its uniqueness under $F_1$
in $M^\ast$. Recall that $D^+$ is $OM$-isomorphic to $M^{\ast\ast}$.
As before, suppose $B^\circ_p=\{\downarrow [p]^n\mid n=0, 1, 2,
\cdots\}$ and $(g, d)$ the fundamental Galois connection between
$I^\sim$ and $M$. Put $v_{[p]}(J)=sup \{n\mid \downarrow
[p]^n\subseteq J\}$ for any $J\in M^\ast$.

\begin{lemma}
Under $F_1$, $v_{[p]}(J)$ exists for any $J\in M^\ast$ and all
$v_{[p]}(J)$ are zero except for a finite number of them. In this
case we have $J=\vee_{[p]\in [P]}\downarrow [p]^{v_{[p]}(J)}$.
\end{lemma}
{\bf Proof.} Assume that $J\in M^\ast$. By (\cite{Deng0}, Lemma 3.9
(1), (2)) $J\subseteq J'\subseteq \downarrow [x]$ for some $J'\in
M^{\ast\ast}, [x]\in G^\ast$ and hence $J\subseteq \downarrow [x]$.
Let $\downarrow [p]^n\subseteq J$. Then $[p]^n\in J$, whence
$[p]^n\leqslant [x]$. Thus we have $v_{[p]}([x])$ as an upper bound
of $\{n\mid \downarrow [p]^n\subseteq J\}$ owing to the fact that
$v_{[p]}([x])$ is the greatest of $n$ such that $[p]^n\leqslant
[x]$. Thus $sup \{n\mid \downarrow [p]^n\subseteq J\}$ exists. And
what is more, all $v_{[p]}([x])$ are zero except for a finite number
of them, so are $v_{[p]}(J)$ because $v_{[p]}(J)\leqslant
v_{[p]}([x])$ for all $[p]\in [P]$. \ \ $\square$

\begin{lemma}
Assume that $F_1$ holds. Then $J=\vee_{[p]\in [P]}\downarrow
[p]^{v_{[p]}(J)}$ with all $v_{[p]}(J)=0$ except for a finite number
of them is unique for any $J\in M^\ast$.
\end{lemma}
{\bf Proof.} Let $[p]\in [P]$. By (\cite{Deng0}, Lemma 3.6)
$\downarrow [p]$ is both a prime element and an atom of
$M^{\ast\ast}$. From (\cite{Deng0}, lemma 3.5) we know that
distributive law of $\cdot$ with respect to $\wedge$ holds. Thus by
Example 2.1 (3) $\downarrow [p]^n\wedge \downarrow [q]^m=\downarrow
[1]$ in $M^{\ast\ast}$ for $[p]\neq [q]$ with $[p], [q]\in [P]$. But
$\downarrow [p]^n\wedge \downarrow [q]^m=\downarrow [1]$ also holds
in $M^\ast$ because by (\cite{Deng0}, Remark 3.1). $i$ in Galois
connection $(i, g\circ d)$ between $M^{\ast\ast}$ and $M^\ast$ is
$\wedge$-preserving. On the other hand, $\downarrow [p]^n$ is
strongly completely $\vee$-irreducible by (\cite{Deng0}, Proposition
2.5 (3)), a fortiori, strongly $\vee$-irreducible. Thus follows
uniqueness of $J=\vee_{[p]\in [P]}\downarrow [p]^{v_{[p]}(J)}$. \ \
$\square$

\begin{prop}
Assume that $F_1$ holds, then each $J\in M^{\ast\ast}$ is a
principal $B$-ideal.
\end{prop}
{\bf Proof.} Let $J\in M^{\ast\ast}$. We have $J=\vee_{[p]\in
[P]}\downarrow [p]^{v_p(J)}$ with all $v_{[p]}(J)=0$ except for a
finite number of them. hence $J=\vee \downarrow
[p_i]^{v_{[p_i]}(J)}$ with $v_{[p_i]}(J)\neq 0$ for $i=1, 2, \cdots,
n$. But $J, \downarrow [p_i]^{v_{p_i}(J)}\in M^{\ast\ast}$ and by
Example 2.2 (1) $\vee \downarrow [p_i]^{v_{p_i}(J)}=\prod \downarrow
[p_i]^{v_{p_i}(J)}=\downarrow \prod [p_i]^{v_{p_i}(J)}$, which is a
principal $B$-ideal. Thus $J$ is a principal $B$-ideal. \ \
$\square$

\begin{prop}
(Supplement to Theorem 2.1) The following conditions are equivalent.

{\rm (1)} $F_1$ holds.

{\rm (2)} each $J\in M^\ast$ can be decomposed into a finite join of
$\downarrow b$ with $b\in B$.

{\rm (3)} $G^\ast$ is $OM$-isomorphic to $\prod^{ed}N_{[p]}$.
\end{prop}
{\bf Proof.} (1) is equivalent to (2). That (1) implies (2) follows
from Lemma 2.5. Conversely by hypothesis for any $a\in G^\ast$ we
have$\downarrow a=\vee \downarrow b_i$ for some $B_a=\{b_i\mid i=1,
2, \cdots, n\}\subseteq B$. Evidently $B_a\in (\downarrow a\cap
B)^{(f)}$. Since $M^{\ast\ast\ast}=\{\downarrow a\mid a\in G^\ast\}$
is order-isomorphic to $G^\ast$ by (\cite{Deng0}, Example 2.1 (2)),
we have $a=\vee b_i$, whence $F_1$ holds.

(1) is equivalent to (3). We know that (1) implies (3) from
Proposition 2.11 (2). Conversely, recall the Galois connection
$(r^\ast_{[p]}, i^\ast_{[p]})$ between $\prod^{ed}N_{[p]}$ and
$N_{[p]}$ (\cite{Deng0}, Definition 3.5).$b_{[p]}\in
\prod^{ed}N_{[p]}$ such that $r_{[q]}(b_{[p]})=1$ for $[q]=[p]$ and
$r_{[q]}(b_{[p]})=0$ for $[q]\neq [p]$ is clearly a prime element of
$\prod^{ed}N_{[p]}$ for any $[p]\in [P]$. For each $a\in
\prod^{ed}N_{[p]}$, we have $a=\vee i^\ast_{[p]}\circ
r^\ast_{[p]}(a)$ (\cite{Deng0}, Lemma 3.9). But $i^\ast_{[p]}\circ
r^\ast_{[p]}(a)=r^\ast_{[p]}(a)b_{[p]}$, whence $a=\vee
r^\ast_{[p]}(a)b_{[p]}$. By the $OM$-isomorphism which transforms
addition of $\prod^{ed}N_{[p]}$ into multiplication of $G^\ast$, (1)
follows. \ \ $\square$

\begin{rem}\rm
Proposition 2.13 (3) can be considered as order representation
theorem of $G^\ast$ under $F_1$. By (\cite{Deng0}, Example 3.3) we
have $M^\ast=\prod^{id}B^\circ_{[p]}$. $\prod^{id}B^\circ_{[p]}$
endowed with the multiplication defined by $a\cdot b=\vee
i_{[p]}(r_{[p]}(a)\cdot r_{[p]}(b))$ for any $a, b\in
\prod^{id}Q^\circ_{[p]}$, becomes an ordered monoid, so does
$M^\ast$. Thus $M^\ast$ as a monoid is monoid-isomorphic to
$\prod^{id}B^\circ_{[p]}$, which is the internal direct product of
monoids $B^\circ_{[p]}$'s in algebraic sense. If we adopt the
notation of join, then $J=(\downarrow [p_1])^{v_{[p_1]}(J)}\vee
\cdots\vee (\downarrow [p_n])^{v_{[p_n]}(J)}$, which means that in
the $T_0$-topological space $X=\{[p]^n\mid [p]\in [P]\}$ (see
\cite{Deng0}, Proposition 2.8) each closed set $J$ which is neither
empty nor the space itself can be written uniquely as a join of a
finite number of point-closures in case of $M$.
\end{rem}

\begin{rem}\rm
In view of Proposition 2.8, 2.9 the main result of
\cite{Martinez}(\cite{Martinez}, Theorem 2.17) is a special case of
equivalence of (6) to (2) in Theorem 2.1, and (\cite{Novak}, Theorem
6.5) is a special case of equivalence of $F_2$ and $D_4$ of (3) to
(2) in Theorem 2.1.
\end{rem}

\section{Krull domains}
 We use $D, D(I^\ast)$ to denote the collection of all divisors, the
 collection of all divisorial ideals associated with $I^\ast$ respectively.

\begin{defn}
By the condition $D_6$ we mean that each $J\in M^{\ast\ast}$ can be
decomposed into a join of a finite number of powers of atoms.
\end{defn}

\begin{lemma}
$R$ is a Krull domain if and only if $D_6$ holds.
\end{lemma}
{\bf Proof.} From (\cite{Bourbaki0}, Chapter VII, $\S$1, n$^\circ$3,
Theorem ) and its proof we know that $R$ is a Krull domain if and
only if each element of $D^+$ can be decomposed into a product of a
finite number of powers of atoms. From (\cite{Bourbaki0}, Chapter
VII, $\S$1, n$^\circ$1, n$^\circ$2) we know that $D^+$ is a lattice
monoid and distributive law of $\cdot$ with respect to $\wedge$
holds. Hence by Example 2.1 (3), Example 2.2 (1), a product can be
reduced to a join and conversely, since $D^+$ is $OM$-isomorphic to
$M^{\ast\ast}$, we can identify $D^+$ with $M^{\ast\ast}$, i.e., we
regard $g(a^-)$ as $div(a^-)$ for each $a^-\in D(I^\ast)$. Thus
$D_6$ holds.\ \ $\square$

~

Now we study decomposition problem in $M^{\ast\ast}$, as lattice
monoid, and put $A_t=$ the set of all atoms of $M^{\ast\ast}$.

From above mentioned we know that in $M^{\ast\ast}$, as a lattice
monoid, distributive law of $\cdot$ with respect to $\wedge$ holds.
By Example 2.3 (1) $e\in A_t$ is a prime element of $M^{\ast\ast}$.
hence $E=\{e^n\mid n=1, 2, \cdots$ and $ e\in A_t\}$ has the same
properties as $B$ in Lemma 2.1. If $D_6$ holds, then as we did in
Proposition 2.11 we can show that $M^{\ast\ast}$ is $OM$-isomorphic
to $\prod^{ed} E_e$, where $E_e=\{e^n\mid n=1, 2, \cdots\}$. We have
the following

\begin{lemma}
If $D_6$ holds, then $M^{\ast\ast}$ is $OM$-isomorphic to
$\prod^{ed}E_e$.
\end{lemma}

\begin{lemma}
Assume that $D_6$ holds. Then

{\rm (1)} $\subseteq$ in $M^{\ast\ast}$ is defined by $\cdot$.

{\rm (2)} Cancellation law and distributive law of $\cdot$ with
respect to $\vee$ hold in $M^{\ast\ast}$.

{\rm (3)} $DCC$ holds in $M^{\ast\ast}$.
\end{lemma}
{\bf Proof.} (1) and (2). As is easily known, $\prod^{ed}E_e$ is
$OM$-isomorphic to $\prod^{ed}N_e$, where $N_e=N$ and in
$\prod^{ed}N_e$ $\leqslant$ is defined by addition and cancellation
law with respect to addition holds. Besides in $N$ distributive law
of $+$ with respect to $\vee$ is $n+max \{m_1, m_2\}=max\{n,
m_1\}+max\{n, m_2\}$, which is trivial. Hence it also holds in
$\prod^{ed}N_e$. Through the $OM$-isomorphism the same is true in
$M^{\ast\ast}$.

(3) $DCC$ holds in $\prod^{ed}N_e$ from the proof of necessity of
(\cite{Bourbaki1}, $\S$1, n$^\circ$13, Theorem 2). $DCC$ also holds
in $M^{\ast\ast}$ by the $OM$-isomorphism. \ \ $\square$

\begin{lemma}
Suppose the conditions (1), (2) and (3) in Lemma 3.3 are all true,
then every $\vee$-irreducible element of $M^{\ast\ast}$ $\in E$.
\end{lemma}
{\bf Proof.} Put $C=\{J\in M^{\ast\ast}\mid J$ is $\vee$-irreducible
but is not an element of $E\}$. Suppose $C\neq \varnothing$. Then by
(3) there exists a minimal $J^\ast$ in $C$. From the proof of
(\cite{Bourbaki1}, $\S$1, n$^\circ$13, Lemma ) it follows that
$J\subseteq J^\ast$ for some $J\in A_t$ and $J\subset J^\ast$
because otherwise $J=J^\ast$ would contradict to $J^\ast\in C$. By
(1) $J^\ast= J\cdot J'$ for some $J'\in M^{\ast\ast}$ and we have
$J'\subset J^\ast$ because $J'=J^\ast$ would lead to $J=\downarrow
[1]$ by cancellation law, a contradiction. Thus $J'\notin C$, whence
either $J'$ is not $\vee$-irreducible or $J'\in E$.

In the former case $J'=J_1\vee J_2$ with $ J_1\subset J', J_2\subset
J'$ and $J^\ast=J\cdot(J_1\vee J_2)=(J\cdot J_1)\vee (J\cdot J_2)$
by the distributive law of $\cdot$ with respect to $\vee$ and we
have $J\cdot J_1\subset J^\ast, J\cdot J_2\subset J^\ast$ by
cancellation law, a contradiction. In the latter case $J^\ast=J'$ or
$J^\ast=J\vee J'$ according as $J'$ is a power of the same atom as
$J$ or not, another contradiction. It follows that $C=\varnothing$.
\ \ $\square$

~

By (\cite{Bourbaki0}, Chapter VII, $\S$1, n$^\circ$1, Proposition 2)
$D^+$ is closed under arbitrary meets and under those joins whose
associated divisorial ideals are not $(0)$, whence $M^{\ast\ast}\cup
\{\downarrow [0]\}$ is a complete lattice.

\begin{prop}
If $R$ is a Krull domain, then $M^{\ast\ast}\cup \{\downarrow [0]\}$
has a topological representation, in which $J\in M^{\ast\ast}$ is a
closed set of some topological $T_0$-space and $e^n$ with $n>0$ its
point-closures.
\end{prop}
{\bf Proof.} By Lemma 3.1 $D_6$ holds in $M^{\ast\ast}$. From
Proposition 2.8 taking $E$ as $B$ we know that $e^n$ is strongly
$\vee$-irreducible. The result follows from (\cite{Deng0}, Lemma
2.1). \ \ $\square$

\begin{thm}
The following conditions are equivalent.

{\rm (1)} $R$ is a Krull domain.

{\rm (2)} $D_6$ holds in $M^{\ast\ast}$.

{\rm (3)} $DCC$ holds in $M^{\ast\ast}$ and each $\vee$-irreducible
element of $M^{\ast\ast}\in E$

{\rm (4)} Distributive law of $\cdot$ with respect to $\vee$ and
cancellation law hold in $M^{\ast\ast}$ and $\subseteq$ is defined
by $\cdot$, and $DCC$ holds in $M^{\ast\ast}$.

{\rm (5)} $M^{\ast\ast}$ is $OM$-isomorphic to $\prod^{ed} N_e$.
\end{thm}
{\bf Proof.} That (1) is equivalent to (2) follows from Lemma 3.1.
By Lemma 3.3, (2) implies (4). By the proof of Lemma 3.3 (2) implies
(5). That (4) implies (2) follows from Lemma 3.4 and Proposition 2.1
(2) taking $E$ as $B$.

(2)is equivalent to (3) This follows from Proposition 2.1(2) because
take $E$ as $B$ and in this case lemma 2.1 (1) (3) taking $E$ as $B$
hold still.

That (5) implies (2) follows from Proposition 2.13 (3) by using
$A_t$ as the index set instead of $[P]$ and noting that $b_e$ is an
atom of $\prod^{ed} N_e$. \ \ $\square$

~

Now we turn to the mutual relations among Krull domains, Dedekind
domains, unique factorization domains and principal ideal domains.

In the following $(g, d)$ will denote the fundamental Galois
connection between $I^\sim$ and $M$.

\begin{lemma}
$D(I^\ast)=I^\ast$ if and only if $d$ is onto.
\end{lemma}
{\bf Proof.} Note that $d(\downarrow [0])=(0)$. Assume that
$D(I^\ast)=I^\ast$ and $a\in I^\ast$. Then we have $g(a)\in
M^{\ast\ast}$ and $d(g(a))=a$, whence $d$ is onto.

Conversely let $a\in I^\ast$. Since $d$ is onto, there exists $J\in
M^\ast$ such that $d(J)=a$. But $d(J)\in D(I^\ast)$ because $J\neq
\downarrow [0]$ and hence $a\in D(I^\ast)$. Thus $I^\ast\subseteq
D(I^\ast)$ and so $I^\ast=D(I^\ast)$. \ \ $\square$

\begin{lemma}
If each integral divisor is principal, then $g$ is onto.
\end{lemma}
{\bf Proof.} Note that $g((0))=\downarrow [0]$. Let $J\in M^\ast$,
then $\cap_{[x]\in J}(x)=(y)$ and hence $\downarrow [y]=\vee_{[x]\in
J}\downarrow [x]=J$  by [1, Example 2.1(1)(2)]. Thus
$g((y))=\downarrow [y]$ and so $g$ is onto. \ \ $\square$

\begin{thm}
{\rm (1)} $R$ is a Dedekind domain if and only if $R$ is a Krull
domain and $D(I^\ast)=I^\ast$.

{\rm (2)} $R$ is a unique factorization domain if and only if $R$ is
a Krull domain and each integral divisor is principal.
\end{thm}
{\bf Proof.} (1) Necessity. Suppose $R$ is a Dedekind domain, then
each $a\in I^\ast$ is inversible  and hence is divisorial by
(\cite{Bourbaki0}, Chapter VII, $\S$1, n$^\circ$2), whence
$D(I^\ast)=I^\ast$. Furthermore, from (\cite{Deng0}, Remark 3.1) we
know that $D(I^\ast)$ ordered by inverse inclusion is
$OM$-isomorphic to $M^{\ast\ast}$ because in this case $D(I^\ast)$
is closed under $\cdot$, whence $D_6$ holds. It follows that $R$ is
a Krull domain by Lemma 3.1.

Sufficiency. By Lemma 3.5 $d$ is onto and so $M^{\ast\ast}$ is
$OM$-isomorphic to $I^\ast$. Since $D_6$ holds in $M^{\ast\ast}$,
$R$ is a Dedekind domain.

(2) necessity. By Lemma 2.5 and its proof , $D_6$ holds , whence $R$
is a Krull domain by Lemma 3.1. Besides, that each $J\in
M^{\ast\ast}$ is a principal $B$-ideal follows Proposition 2.12.

Sufficiency. In this case $M^{\ast\ast}$ is $OM$-isomorphic to
$G^\ast$, whence $D_6$ in $M^{\ast\ast}$ is reduced to $F_1$ in
$G^\ast$. Hence $R$ is a unique factorization domain. \ \ $\square$

\begin{cor}
$R$ is a principal ideal domain if and only if $R$ is both a
Dedekind domain and a unique factorization domain.
\end{cor}
{\bf Proof.} Sufficiency. By Theorem 3.1 (1)(2) each ideal is
divisorial, which is principal. Thus $R$ is a principal ideal
domain.

Necessity. We have $I^\ast=D(I^\ast)=L^\ast$ and $R$ is a unique
factorization domain (see Remark 3.1 below), whence $R$ is a Krull
domain by Theorem 3.2 (2), which together with $D(I^\ast)=I^\ast$
implies that $R$ is a Dedekind domain by Theorem 3.2(1). \ \
$\square$

\begin{rem}\rm
Any principal ideal domain is a unique factorization domain. In
fact, as is well known, $DCC$ holds in $I^\ast$ and $I^\ast$ is a
lattice, so $DCC$ holds in $G^\ast$ and $G^\ast$ is a lattice
because $I^{\ast \sim}$ is $OM$-isomorphic to $G^\ast$. Besides,
$G^\ast$ is the integral part of lattice group $K^\ast/U$
(\cite{Bourbaki1}, $\S$1, n$^\circ$9, Proposition 8) in notation of
multiplication and so the conditions (1)(2) of Lemma 3.3 are
satisfied, whence each $\vee$-irrducible element of $G^\ast$ belongs
to $B$ by Lemma 3.4 applied to $G^\ast$ with $B$ instead of $Q$. It
follows that $F_1$ holds by Proposition 2.1 (2).
\end{rem}

Note that in case of $R$ being a Dedekind domain or a unique
factorization domain or a principal ideal domain, $g$ is injective
or surjective or bijective.

\begin{rem}\rm
We have already known from Remark 2.2 that Proposition 2.13 (3) can
be viewed as the order representation of $G^\ast$ under $F_1$.

(1)(topological representation of $G=G^\ast\cup \{[0]\}$ under
$F_1$) Assume that $F_1$ holds. By Theorem 3.2 (2) each integral
divisor is principal, whence $g$ is onto by Lemma 3.6, whose proof
shows each $J\in M^\ast$ is a principal $B$-ideal. Thus by
(\cite{Deng0}, Proposition 2.7), which is true under $F_1$ because
$F_1$ implies $D_1$, $B_2$ holds. We claim that $G$ is a complete
lattice. In fact, let $\{a_t\mid t\in T\}\subseteq 2^{G^\ast}$. We
have $a_t=\vee A_t$ for some $A_t\in (\downarrow a_t\cap B)^{(f)}$
by $F_1$ and $\vee \cup A_t$ exists by $B_2$ so that $\vee a_t$
exists and is equal to $\vee \cup A_t$. Thus $G$ is a complete
lattice. By Proposition 2.8 each $b\in B$ is strongly
$\vee$-irreducible and hence by (\cite{Deng0}, Lemma 2.1) $G$ has a
topological representation.

(2)(algebraic interpretation of $G^\ast$ under $F_1$) By proposition
2.10 $G^\ast=\prod^{id}Q_{[p]}$, which is $OM$-isomorphic to
$\prod^{ed}Q_{[p]}$. For the Latter $G^\ast$, as a monoid, is the
direct product of monoids $Q_{[p]}$'s in a algebraic sense. For the
former if notation of addition is adopted, for each $[x]\in G^\ast$
we have $[x]=v_{[p_1]}([x])\cdot [p_1]+\cdots +v_{[p_n]}([x])\cdot
[p_n]$, whence $G^\ast$ is an analogue of a module for the
background of $N$ with $[P]$ as a base. If we adopt notation of
join, then $[x]=[p_1]^{v_{[p_1]}([x])}\vee \cdots \vee
[p_n]^{v_{[p_n]}([x])}$, which according to (1) means any closed set
$[x]$, which is neither empty nor the space itself, can be written
uniquely as a join of a finite number of point-closures.

(3)(algebraic interpretation of a subposet of the first kind or of
the second kind under $F_1$). Note first that $Q_{[p]}$ is a
$B$-ideal, which is a subposet of the first kind of $G^\ast$, as can
be easily verified. Now consider $G^\ast=\prod^{id}Q_{[p]}$. We have
$G^\ast=\downarrow [0]\in M$ and $\prod^{id}Q_{[p]}=\vee (\cup
Q_{[p]})=\vee Q_{[p]}$ (see \cite{Deng0}, remarks below Proposition
2.2). Evidently $Q_{[p]}\cap Q_{[q]}=\downarrow [1]$ for $[p]\neq
[q]$.

Put $P_{[q]}=\vee_{[p]\neq [q]}Q_{[p]}$ and we claim that $\cap
P_{[q]}=\downarrow [1], P_{[p]}\vee P_{[q]}=\downarrow [0]$ for
$[p]\neq [q]$ and $P_[q]$ is of the second kind. In fact, let
$A_{[q]}=\{Q_{[p]}\mid [p]\neq [q]\}$ for $[q]\in [P]$. We have
$\cap P_{[q]}=\cap \vee A_{[q]}$, which by complete distributivity
(\cite{Deng0}, Proposition 2.7) is equal to $\vee \{\cap f([P])\mid
f\in \prod A_{[q]}\}=\vee \{[1]\}=\downarrow [1]$. And $P_{[p]}\vee
P_{[q]}=\vee Q_{[p]}=G^\ast=\downarrow [0]$ for $[p]\neq [q]$. At
last, since $P_{[q]}$ is closed under $\wedge$, it is of the second
kind.

On the other hand, it will be easily verified that
$P_{[q]}=J_{[q]}$, whence $J_{[q]}$ is of the second kind. We have
$J_{[q]}=\{[x]\mid v_{[q]}([x])=0\}$, which by (\cite{Deng0},
Proposition 2.1(1)) is a prime $B$-ideal because $G^\ast\setminus
J_{[q]}$ is clearly a $B$-filter.
\end{rem}

\begin{rem}\rm
(order or algebraic or topological interpretation of $M^{\ast\ast}$
under $D_6$) By Lemma 3.2 and remarks above it,
$M^{\ast\ast}=\prod^{ed} E_e$, which is $OM$-isomorphic to
$\prod^{ed} E_e$, which in turn is $OM$-isomorphic to $\prod^{ed}
N_e$, where $N_e=\{0, 1, 2, \cdots\}$ ordered by the order defined
by addition. $M^{\ast\ast}$, as a poset, is order-isomorphic to
$\prod^{ed}N_e$, which can be viewed as the order representation of
$M^{\ast\ast}$. Besides, $M^{\ast\ast}$, as a monoid, is
monoid-isomorphic to $\prod^{ed} E_e$ and so $M^{\ast\ast}$ is the
direct product of monoids $E_e$'s in algebraic sense. If we adopt
notation of addition in case of $M^{\ast\ast}=\prod^{id}E_e$, then
for each $J\in M^{\ast\ast}$ we have $J=v_{e_1}(J)\cdot e_1+\cdots
+v_{e_n}(J)\cdot e_n$, whence $M^{\ast\ast}$ is an analogue of a
module for the background of $N$ with $A_t$ as a base. Finally by
Proposition 3.1 $M^{\ast\ast}\cup \{\downarrow [0]\}$ has a
topological representation with $E$ as the space $X$ and we have
$J=e_1^{v_{e_1}(J)}\vee \cdots \vee e_n^{v_{e_n}(J)}$, which means
any closed set $J$ of $X$ which is neither empty nor the space
itself can be written uniquely as a join of a finite number of
point-closures.
\end{rem}

\begin{rem}\rm
Each element of $L^\ast$ is regular in \cite{Bourbaki2}. By
symmetrization (\cite{Bourbaki2}, $\S$1, n$^\circ$4, Theorem 1) we
extend $I^\ast$ to the collection $F$ of fractional ideals, in which
each element of $L^\ast$ is inversible. Completely integral
closeness means each divisor being inversible (\cite{Bourbaki0},
Chapter VII, $\S$1, n$\circ$2, Theorem 1), while a Dedekind domain
is equivalent to each element of $F$ being inversible. Thus in Lemma
3.3 cancellation law can deduce (1) and the remaining part of (2).
\end{rem}

\begin{center}
College of Mathematics and Econometrics, \\Hunan University,
Changsha 410082, Hunan, China.
\end{center}